%% file: TSTsArXiv.tex
\pdfoutput=1
\documentclass[reqno]{amsart}
   \newif\ifTeXpad
   \TeXpadfalse

\usepackage{amssymb,mathrsfs}
\usepackage{graphicx,xcolor}
\usepackage[bookmarks,hyperfootnotes=false]{hyperref}
\usepackage[nameinlink,capitalise,noabbrev]{cleveref}
\usepackage[msc-links]{amsrefs} 
\usepackage{doi}
   \input letterswitharrowsRhd
\usepackage{enumitem}
\usepackage{comment}

\colorlet{darkishRed}{red!80!black}
\colorlet{darkishBlue}{blue!60!black}
\colorlet{darkishGreen}{green!60!black}
\hypersetup{
  draft=false,
  bookmarksopen=true,
  colorlinks,
  linkcolor=darkishRed,
  citecolor=darkishGreen,
  urlcolor=darkishBlue
}

\renewcommand{\PrintDOI}[1]{\doi{#1}}

\renewcommand{\leq}{\leqslant}
\renewcommand{\geq}{\geqslant}

\newcommand{\R}{\mathbb{R}}
\newcommand{\N}{\mathbb{N}}
\newcommand{\B}{\mathcal{B}}
\newcommand{\C}{\mathcal{C}}
\newcommand{\F}{\mathcal{F}}
\renewcommand{\P}{\mathcal{P}}
\newcommand{\RR}{\mathcal{R}}
\newcommand{\T}{\mathcal{T}}
\newcommand{\td}{tree-decom\-po\-si\-tion}
\newcommand{\FBk}{\B_k}
\newcommand{\FCn}{\C_n}
\newcommand{\FP}{\P_{\rm s}}

\newtheorem{theorem}{Theorem}[section] 
\newtheorem{corollary}[theorem]{Corollary}
\newtheorem{lemma}[theorem]{Lemma}

\theoremstyle{definition}
\newtheorem{definition}[theorem]{Definition}

\setlist[enumerate]{
  label=\textrm{(\roman*)},   
  font=\normalfont            
}
\def\COMMENT{}
\def\ucl(#1){\lfloor #1 \rfloor}

\title{Tangle structure trees}
\author{Hanno von Bergen and Reinhard Diestel}
\date{\today}

\begin{document}

\maketitle

\begin{abstract}
We introduce a comprehensive data structure, tangle structure trees, which simultaneously displays all the $\F$-tangles of an abstract separation system for very general obstruction sets~$\F$. It simultaneously also displays certificates $\sigma\in\F$ for any non-existence of such tangles, or for the non-extend\-ability of low-order tangles to higher-order ones.

Our theorem can be applied to produce the structures%
   \COMMENT{}
   of the classical tree-of-tangles and tangle-tree duality theorems, both for graph tangles and for their known generalizations to more general separation systems. It extends those theorems to obstruction sets~$\F$ that need not define profiles (as they must in all known trees of tangles) or consist of stars of separations (as they must in traditional tangle-tree duality).

Our existence proof for these structure trees is constructive. The construction has been implemented in open-source software available for tangle detection and further analysis.
\end{abstract}

\medskip\section{Introduction}

\noindent
   The notion of `tangles' was originally introduced by Robertson and Seymour~\cite{GMX} as an abstract concept of high local connectivity in graphs, one that unifies several more concrete such notions, such as highly connected subgraphs or minors. What all these notions have in common is that, given any low-order separation of the graph, any highly cohesive substructure must lie mostly on one of its two sides: since the separation has low order, it cannot split it into two roughly equal halves. A~tangle remembers only how all the low-order separations are oriented in this way, each towards that highly cohesive substructure: the collection of all these oriented separations is called a {\em tangle\/}. See~\cite{DiestelBook25} for a precise definition and basic facts about graph tangles.

Tangles have since been generalized to more general settings than graphs. In the setting of set partitions,%
   \COMMENT{}
    they offer a precise theoretical basis for `fuzzy' real-world problems such as clustering in large datasets~\cite{TangleBook}. All these, including graphs, are special cases of so-called {\em abstract separation systems\/}~\cite{ASS}: an axiomatic setting that assumes only the most basic properties of graph separations. This is the setting in which tangles are most easily and comprehensively treated, and we shall use this framework also in this paper.

There are two fundamental theorems about tangles, which both have their origins in~\cite{GMX} and are also treated in~\cite{DiestelBook25}. The {\em tree-of-tangles\/} theorem exposes a tree-like structure in the graph or dataset whose tangles we consider, which `distinguishes' the tangles in that they are shown to `live in' different areas of this tree-like structure. The other is the {\em tangle-tree duality\/} theorem, which certifies the non-existence of a tangle by exposing another tree-like structure in the graph or dataset, one whose existence clearly precludes the existence of a tangle.%
   \footnote{In a graph, this would be a tree-decomposition into parts too small for a tangle to live in.}
   This is the more fundamental of the two theorems, in the sense that the other can be reduced to its abstract version (see below) but not vice versa~\cite{ToTfromTTD} .

Both the above theorems have been generalized to  $\F$-tangles of abstract separation systems, ways of simultaneously orienting the `separations' of a given structure in such a way that no three of those oriented separations form an element of~$\F$.%
   \footnote{In the case of graph tangles, $\F$~consists of the sets of up to three oriented separations of the graph such that the union of the `back' sides of these three separations covers the entire graph~\cite{DiestelBook25}.}
   This collection~$\F$ has to satisfy some constraints that depend on which of the two theorems we are considering.

For the (known) tree-of-tangles theorems, those which display all the $\F$-tangles of a separation system simultaneously, $\F$~has to include triples of oriented separations that are reminiscent of ultrafilters. If our separations are bipartitions of a fixed set~$V\!$, for example, then $\F$ must contain all triples of the form~$\{A,B,\overline{A\cap B}\}$: if $A$ and $B$ are subsets of~$V\!$ deemed `large', in that the $\F$-tangle orients the partitions $\{A,\overline A\}$ and $\{B,\overline B\}$ towards them, then the complement of $A\cap B$ cannot also be `large'. Such $\F$-tangles are called {\em profiles\/}, and all existing tree-of-tangle theorems for $\F$-tangles assume that these are profiles~\cite{ProfilesNew,CDHH13CanonicalAlg,CDHH13CanonicalParts,carmesin2022entanglements,JoshRefining,SARefiningInessParts,SARefiningEssParts,AbstractTangles}. There are even a tree-of-tangles theorems that display profiles of different `order' simultaneously. Such profiles, however, must be {\em robust\/}: they also must not contain certain triples similar to, but slightly different from, those above. More about profiles in \cref{sec:apps}.%
   \footnote{In a follow-up paper~\cite{TSTs2} we show how our structure trees give rise to trees of $\F$-tangles whenever these are robust, even if they are not profiles.}

For the tangle-tree duality theorems, $\F$~has to satisfy another constraint: its elements have to be {\em stars\/}, nested sets of oriented separations pointing towards each other~\cite{ASS}. The tree-like structures by which tangle-tree duality theorems certify the non-existence of $\F$-tangles require that $\F$ consist of stars: if it does not, these structures cannot be tree-like.

This state of the art leaves $\F$-tangles that are not profiles without any known way of displaying all those tangles simultaneously. And if $\F$ does not consist of stars, it leaves separation systems that have no $\F$-tangles without any known way of organizing the elements of~$\F$ into a data structure that displays them as easily checkable certificates for the non-existence of such tangles.

The structure trees whose existence we prove in this paper do both these things: they display all the $\F$-tangles of a separation system even when they are not profiles, and for those orientations of the separation system that are not $\F$-tangles they display certificates in~$\F$ that show why they are not.%
   \footnote{In a follow-up paper~\cite{TSTs2} we show that if $\F$ does consist of stars, our duality theorem implies the same tangle-tree dichotomy as the classical one.}

Tangle structure trees display all this information in a single, comprehensive, data structure, which is maximally efficient in the following, structural, sense. A~single tangle is most efficiently displayed by listing just those of its elements that are minimal in the poset of oriented separations that comes with a separation system. Indeed, all these are needed to determine the tangle, but any other elements can be deduced from them. Our structure trees display all the $\F$-tangles of a separation system~$S$ simultaneously by listing only the oriented separations that are minimal in one of those tangles, including the tangles of subsystems consisting of separations of lower order. And for orientations of~$S$ that are not $\F$-tangles it displays certificates from~$\F$ inside those same minimal subsets of oriented separations,%
   \COMMENT{}
   in a way readily accessible in the structure tree.

Our paper is organized as follows. In \cref{sec:defs} we  provide the necessary background for tangles in abstract separation systems, the framework in which we shall construct our structure trees and prove their existence.

In~\cref{sec:TSTs} we introduce basic tangle structure trees. In~\cref{sec:TSTexistence} we prove their existence, and in \cref{sec:efficient} we show how to make them efficient. In \cref{sec:MainThm} we collect all this information together to prove our main results, the `$\F$-tangle structure theorems'. In \cref{thm:generalduality2} we note specifically how the main structure theorem provides certificates from~$\F$ for the non-existence of $\F$-tangles if there are none.\looseness=-1

In \cref{sec:apps}, finally, we apply our results to three particular types of $\F$-tangles. The first two include $\F$-tangles in graphs: those induced by $k$-blocks, and those that are profiles. In neither of these does the defining set~$\F$ consist of stars, so traditional tangle-tree duality fails. We derive new dichotomy theorems for these.

In our third application we treat $\F$-tangles that encode clusters in large datasets. These are neither profiles nor do their obstruction sets~$\F$ consist of stars. Our structure trees display all the clusters by way of their $\F$-tangles, while also witnessing the absence of clusters from the other areas of the data set by displaying elements of~$\F$.\looseness=-1

\medskip\section{Tangle basics}\label{sec:defs}

\noindent
   In this section we give precise definitions and notation for abstract tangles, largely following \cite{ASS} and indicating any deviations.%
   \footnote{The most important difference is that, for historical reasons, the partial ordering on~$\vS$ used in~\cite{ASS} is the inverse of ours. So terms like `large' and `small', infima and suprema etc, are reversed.}

Tangles of graphs are ways of orienting their separations, each towards one of its two sides. Abstract tangles are designed to work in scenarios where there need not be anything to `separate'. In order to retain our intuition from graphs, however, we continue to refer to the things of which our abstract tangles pick one of two variants (which they will indeed do) as `separations'. These are defined by noting some key properties of graph separations and making them into axioms, as follows.

A \emph{separation system} $(\vS, \leq, ^*)$ is a set~$\vS$, whose elements we call \emph{oriented separations}, that comes with a partial ordering~$\leq$ on~$\vS$ and an order-reversing involution $^*\colon \vS \to \vS$. Thus, for any two elements%
   \footnote{We often denote the elements of~$\vS$ by letters with an arrow, in either direction, precisely in order to have a simple way to refer to their dual elements: by reversing the arrow. But the arrow directions have no meaning: an arbitrary element of~$\vS$ could be denoted equally as $\vs$ or as~$\sv$.}
   $\vr,\vs$ of~$\vS$ with $\vr \leq \vs$ we have $\vr{}^* \geq \vs{}^*$. We write $\vs{}^*=:\sv$, and call $\sv$ the \emph{inverse} of $\vs$. While we allow formally that $\vs = \sv$, in which case we call $\vs$ and~$s$ \emph{degenerate}, this does not happen often in practice.%
   \footnote{The only degenerate separation of a graph $G=(V,E)$, for example, is~$\{V,V\}$.}
 
If a separation system~$\vU$ happens to be a lattice, that is, if there is a supremum $\vr\lor\vs$ and an infimum $\vr\land\vs$ in~$\vU$ for every two elements $\vr,\vs\in\vU$, we call~$\vU$ a {\em universe\/} of separations. It is {\em distributive\/} if it is distributive as a lattice. A~separation system $\vS\subseteq\vU$ is {\em submodular\/} if for every two elements of~$\vS$ either their infimum or their supremum in~$\vU$ also lies in~$\vS$.

Very rarely we may have separations $\vs\le\sv$; then $\vs$ is {\em small\/} and $\sv$ is {\em large\/}.%
   \footnote{The small separations of a graph~$G$ are those of the form~$(V,A)$ with $A\subseteq V = V(G)$.}
   We say that $\vs$ is \emph{trivial} (and $\sv$ is \emph{co-trivial}) in $\vS$ if there exists a pair of inverse separations $\vr,\rv<\vs$ in~$\vS$. Trivial separations are clearly large, so co-trivial ones are small, but the converse need not hold. See~\cite{ASS} for more on these technicalities if desired.

The set of {\em unoriented separations} in~$(\vS,\le,^*)$ is
$$S:=\{\{\vs,\sv\}:\vs\in\vS\}.$$
We call the elements $\vs,\sv$ of~$s$ its \emph{orientations}. An \emph{orientation of~$S$} is a set $\tau\subseteq\vS$ that contains exactly one orientation of every $s\in S$. For $s\in S$ we then denote by $\tau(s)$ the unique orientation of $s$ contained in~$\tau$.%
   \COMMENT{}
   An~orientation of a subset of~$S$ is a {\em partial orientation\/} of~$S$.

If $\vr\geq\vs$ we say that $\vr$ \emph{points towards} $s$ (and that $\rv$ \emph{points away from} $s$). We say that $\vr$ \emph{points towards} an oriented separation $\vs$ whenever it points towards $s$, i.e., if $\vr \geq \vs$ or $\vr \geq \sv$, and similarly for `points away from'. 
A \emph{star\/} is a set~$\sigma$ of non-degenerate oriented separations that point towards each other. As is easy to check, this happens if and only if $\vr\geq\sv$ (and hence $\vs\geq\rv$) for all  distinct $\vr,\vs\in\sigma$.\looseness=-1

\begin{figure}[ht]
 \center\vskip-6pt
   \includegraphics[scale=1]{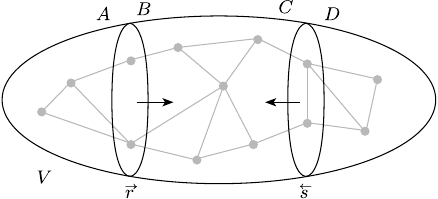}
\vskip-6pt
\caption{\small Nested separations $r = \{A,B\}$ and $s = \{C,D\}$ of a graph. Their orientations $\vr = (A,B)$ and $\sv = (D,C)$ point towards each other, since $\vr\ge \vs$ (as $B\supseteq D$) and $\sv\ge\rv$ (as $C\supseteq A$).}
\label{fig:nested}\vskip-3pt
\end{figure}

Two separations $s,r\in S$ are \emph{nested} if they have orientations that are comparable under $\leq$. Oriented separations are \emph{nested} if their underlying unoriented separations are nested. A subset of~$\vS$ is \emph{nested} if its elements are pairwise nested.

A subset of~$\vS$ is \emph{consistent} if no pair of its elements $\vr,\vs$ with $r \neq s$ point away from each other. Stars are examples of consistent nested sets of oriented~separations.

We shall often be interested in consistent orientations of~$S$. For each of its elements~$\vr\!$, a~consistent orientation of~$S$ will also contain every $\vs > \vr$ other than, possibly, $\vs=\rv$.%
   \COMMENT{}
   Consistent partial orientations of~$S$ are easily seen to extend to consistent orientations of~$S$, unless they contain a separation that is co-trivial in~$\vS$; see~\cite{ASS}.%
   \COMMENT{}

If $\sigma\subseteq\vS$ is consistent, we say that $\vs\in\vS$ is \emph{required} by $\sigma$ if $\vs \notin \sigma$ and $\sigma\cup\{\sv\}$ is inconsistent. We shall see in \cref{lem:closure-consistent} that, pathological cases aside, $\sigma \cup \{\vs \}$ will then be consistent. The \emph{closure} of $\sigma$ is
$$
\lfloor\sigma\rfloor := \sigma \cup\{\vs\in\vS:\text{\(\vs\) is required by \(\sigma\)}\}.
$$
Note that $\sigma$ requires $\vs \notin \sigma$ if and only if there exists an $\vr \in \sigma$ such that $r \neq s$ and $\vs > \vr$. Thus
 $$\lfloor \sigma \rfloor = \sigma \cup \{\, \vs \in \vS: \exists \vr \in \sigma \text{ such that } r \neq s \text{ and }\vs > \vr\}\text{,}$$
 which motivates the term of (upward) `closure'. If $\ucl(\sigma)$ is consistent, then $\lfloor \lfloor \sigma \rfloor \rfloor$ is (defined and) easily shown to be equal to~$\ucl(\sigma)$.%
   \COMMENT{}

If $\sigma$ contains no small separations, the expression above simplifies to
 $$\ucl(\sigma) = \{\,\vs\in\vS : \exists \vr\in\sigma\text{ such that }\vr\le\vs\}.$$
   Indeed, any~$\vs$ in this latter set either lies in~$\sigma$ or there is some $\vr\in\sigma$ such that $\vr < \vs$; in that case $r\ne s$, since otherwise $\sv = \vr < \vs$, making $\vr\in\sigma$ is small.

\begin{lemma}\label{lem:closure}
Let $\sigma\subseteq\tau\subseteq\vS$ be consistent sets. Then $\lfloor\sigma\rfloor\subseteq\lfloor\tau\rfloor$. If $\tau$ is an orientation of all of~$S$, then $\ucl(\tau) = \tau$.
\end{lemma}

\begin{proof}
If $\vs$ is required by $\sigma$, then $\sigma\cup\{\sv\}$ is inconsistent. Since $\sigma\subseteq\tau$, it follows that $\tau\cup\{\sv\}$ is also inconsistent; hence $\vs$ is either an element of or required by~$\tau$. Thus every element of $\lfloor\sigma\rfloor$ lies in $\lfloor\tau\rfloor$, proving the first claim.

If $\tau$ is an orientation of all of~$S$ then $\lfloor\tau\rfloor=\tau$, since for any $\vs\in\ucl(\tau)\setminus\tau$ the set $\tau\cup\{\sv\} = \tau$ would be inconsistent.%
   \COMMENT{}
\end{proof}

Consistent orientations of~$S$ can contain small separations. But examples are rare and can be counter-intuitive, so we often exclude them.%
   \COMMENT{}
   Co-trivial separations cannot occur in consistent orientations of~$S$. Indeed if $\sv$ is co-trivial, witnessed by~$r\in S$, then every orientation of~$S$ will have to orient~$r$ too. But it cannot do so consistently with~$\sv$, since both $\vr$ and~$\rv$ are inconsistent with~$\sv$.%
   \COMMENT{}
    Similarly, $\ucl(\{\sv\})$~is then inconsistent, since it contains $\vr$ and~$\rv$, which are both inconsistent with~$\sv$.

\begin{lemma}\label{lem:closure-consistent}
Let $\tau\subseteq\vS$ be consistent, and suppose that $\tau$ has no elements that are co-trivial in $\vS$. Then $\lfloor\tau\rfloor$ is consistent, and $\ucl(\tau)\setminus\tau$ contains at most one orientation of any $s\in S$.%
   \COMMENT{}
\end{lemma}

\begin{proof}
   If $\lfloor\tau\rfloor$ is inconsistent, it contains orientations $\vr, \vs$ of distinct $r, s \in S$ with $\vs < \rv$. Since $\tau$ is consistent, at least one of $\vr,\vs$ is not in $\tau$. Without loss of generality assume $\vr\notin\tau$: Then $\vr \in \lfloor \tau \rfloor \setminus \tau$, so $\vr$ is required by $\tau$ (by definition of~$\lfloor \tau \rfloor$), which means that $\tau\cup\{\rv\}$ is inconsistent. Thus there exists $\vt\in\tau$ with $t \neq r$ and $\rv<\tv$. We thus have $\vs < \rv <\tv$, that is $\vs$ and $\vt$ point away from each other.\looseness=-1

If $\vs = \vt$, the above inequality yields $\vr, \rv < \tv$, which makes $\vt \in \tau$ co-trivial, contradicting our assumptions. Hence $\vs \neq \vt$ as well as $\vs \neq \tv$, giving $s \neq t$. The fact that $\vs$ and $\vt \in \tau$ point away from each other thus implies that $\{\vs, \vt\}$ is inconsistent, so $\vs \notin \tau$ as well as $\vr\notin\tau$. 

As earlier with $\vr$, the fact that $\vs \in \lfloor \tau \rfloor$ now implies that $\vr < \sv < \tvdash$ for some $\vtdash \in \tau$. The two inequalities together now give $\vtdash < \rv < \tv$, so $\vtdash$ and $\vt$ point away from each other. As $\tau$ is consistent and contains both, this means that $t' = t$. As $\vtdash < \tv$, we thus have $\vtdash = \vt$. But now $\vr < \tvdash = \tv$ as well as $\rv < \tv$. This makes $\vt \in \tau$ co-trivial in~$\vS$, contradicting our assumptions. This proves the first assertion.

For the second, note that~$\tau$ requires any $\vs\in\ucl(\tau)\setminus\tau$, by definition of~$\ucl(\tau)$. This means that $\tau\cup\{\sv\}$ is inconsistent, so $\sv\notin\ucl(\tau)$ since $\ucl(\tau)$ is consistent.
\end{proof}

An \emph{order function} on $S$ is any map $S\to\mathbb{R}$. Unless otherwise mentioned, we denote such order functions as $s\mapsto |s|$. We extend them to~$\vS$ by letting $|\vs| := |\sv| := |s|$. A~separation system~$(\vS,\le,^*)$ given with an order function on~$S$ is an {\em ordered separation system\/}.%
   \COMMENT{}
   For every $k\in\mathbb{R}$ we let
  $$\vSk := \{\vs\in\vS: |s|<k\};$$
this is again an ordered separation system. 

We are sometimes interested in orientations of~$S$ that do not have certain subsets. We typically collect those together in some set~$\F$, whose elements we call \emph{forbidden subsets\/}. Formally, if $\F$~is any set,%
   \COMMENT{}
   we say that $\tau\subseteq\vS$ {\em avoids\/}~$\F$ if $\tau$ has no subset in~$\F$, i.e., if no subset of~$\tau$ is an element of ~$\F$.

\begin{definition}\label{def:tangle}
An $\F$\emph{-tangle} of~$S$ is an $\F$-avoiding consistent orientation of~$S$. The $\F$-tangles of the subsets~$S_k$ of~$S$ are the $\F$-tangles \emph{in}~$\vS$.
\end{definition}

\noindent
   When the context is clear, e.g.\ when $\F$ is some given set as above, we shall usually say `tangle' rather than `$\F$-tangle' in this paper. A~separation $s\in S$ {\em distinguishes\/} two tangles if they orient it differently: if one contains~$\vs$, the other~$\sv$.

\medskip\section{Tangle structure trees}\label{sec:TSTs}

\noindent
   We adopt the graph-theoretic terminology of \cite{DiestelBook25}. A \emph{tree} is a connected acyclic graph. Given nodes~$t,t'$ of a tree $T$ we write $tTt'$ for the unique $t$--$t'$ path in~$T$. A~\emph{rooted tree} is a tree with a distinguished node called its \emph{root}. Given a tree~$T$ with root $r$, we define a partial order $\leq_r$ on $V(T)$ by declaring $x \leq_r y$ if $x$ lies on the path in~$T$ from $r$ to~$y$. Maximal elements, including the root if $|T|=1$, are called \emph{leaves}. Any direct successors in~$<_r$ of a node of~$T$ are its \emph{children}. We write $E_v$ for the set of edges of~$T$ from a node~$v$ to its children.

Let $(\vS,\le,^*)$ be a separation system, and let $\F$ be any set.

\begin{definition}\label{def:septree}
    A \emph{separation tree} $(T,r,\beta)$ on~$\vS$%
   \COMMENT{}
   consists of a rooted tree $(T,r)$ together with an edge labelling $\beta\colon E(T)\to\vS$ such that for every non-leaf ${v\in V(T)}$ there exists a separation $s_v\in S$ such that $\beta$ restricts to a bijection $E_v\to\{\vsv,\svv\}$ and $s_u\neq s_v$ whenever $u <_r v$.
\end{definition}

Thus, every non-leaf node $v$ of a separation tree has either one or two children.  If there is no need to refer to $r$ or $\beta$ explicitly, we usually abbreviate $(T, r, \beta)$ to~$T\!$. For every node~$v$, we write~$\beta_v$ for the set~$\beta\big(E(rTv)\big)$ of edge labels on the path~$rTv$.

\begin{lemma}\label{lem:orichoosesleaf}
    Let $(T,r,\beta)$ be a separation tree on~$\vS$. Every orientation $\tau$ of~$S$ contains $\beta_{\ell}$ as a subset for some unique leaf $\ell\in V(T)$.
\end{lemma}

\begin{proof}
Since $\beta_r=\emptyset\subseteq\tau$, there exists a node $v\in V(T)$ that is maximal with respect to $\leq_r$ subject to $\beta_v\subseteq\tau$. If $v$ is not a leaf, let $v'$ be its child with $\beta(vv')=\tau(s_v)=:\vs$. As $\beta_v\subseteq\tau$, and $\vs\in\tau \setminus \beta_v$ since $s_u\ne s_v$ for all $u<v$, we see that $\beta_{v'}=\beta_v\cup\{\vs\}\subseteq\tau$ contradicts the maximality of~$v$. Hence $v =:\ell$ is our desired leaf.

$\!$For a proof that $\ell$ is unique, let $\ell'$ be another leaf and let $v$ be the greatest common ancestor of $\ell$ and~$\ell'$. Then $v$ has distinct children $w$ and $w'$ with ${\beta(vw)\in\beta_{\ell}}$ and $\beta(vw')\in\beta_{\ell'}$. Since $\beta(vw) \neq \beta(vw')$, they cannot both lie in $\tau$. As $\beta(vw)\in\beta_{\ell}\subseteq\tau$ we must have $\beta(vw')\notin\tau$, and hence $\beta_{\ell'}\not\subseteq\tau$.
\end{proof} 

A separation tree $T$ is \emph{consistent} if $\beta_v\subseteq\vS$ is consistent for every node $v\in T$.

\begin{lemma}\label{lem:notyetnamedlemma}
    Let $(T,r,\beta)$ be a consistent separation tree on~$\vS$, and let $v$ be a non-leaf node of~$T$. Then $\vsv$ and $\svv$ are $\le$-minimal in $\beta_v \cup \{\vsv\}$ and in $\beta_v \cup \{\svv\}$, respectively.
\end{lemma}

\begin{proof}
Suppose not. Without loss of generality assume there exists $\vs\in\beta_v$ with $\vs < \vsv$. Let $w$ be the child of $v$ with $\beta(vw)=\svv$. Then $\{\vs, \svv\} \subseteq \beta_w$ is inconsistent,%
   \COMMENT{}
   contradicting the consistency of $(T,r,\beta)$.
\end{proof}

\begin{corollary}\label{cor:B}
    For every consistent separation tree and non-leaf node~$v$ we have $\{\vsv, \svv\} \cap \lfloor \beta_v \rfloor = \emptyset$. That is, $s_v$~is not oriented by $\lfloor \beta_v \rfloor$.
\end{corollary}

\begin{proof}
Suppose, say, that $\vsv\in\ucl(\beta_v)$. Since $s_u\ne s_v$ for all nodes~$u<v$, we know that $\vsv\notin\beta_v$. Our assumption that $\vsv\in\ucl(\beta_v)$ thus means that $\vsv$ is required by~$\beta_v$. Hence there exists $\vr\in\beta_v$ such that $\vr < \vsv$ (and $r\ne s_v$). This contradicts \cref{lem:notyetnamedlemma}.\looseness=-1
\end{proof}

\cref{lem:notyetnamedlemma} and \cref{cor:B} imply that consistent separation trees offer an efficient data structure for storing consistent orientations of~$S$, in particular, for tangles. Indeed, any such orientation $\tau$ labels the tree's edges between the root and some leaf $\ell$, but this $\beta_{\ell}$ is only a small subset of $\tau$: elements $\vs$ of $\tau$ required by~$\beta_v$ for some $v < \ell$ will not appear in $\beta_{\ell} \setminus \beta_{v}$. But they can be reconstructed from~$\beta_v$ by the consistency of~$\tau$. 

By definition of separation trees, the separations labelling their edges do not repeat along any root-to-leaf path: if $u<v$ then $s_u\neq s_v$. Hence the length of any path from the root to a leaf is at most $|S|$. Since every node has at most two children, this bounds the size of the tree:

\begin{corollary}\label{lem:smoll}
 Separation trees on~$\vS$ have at most $2^{|S|}$ leaves and fewer than $2^{|S| + 1}$ nodes. \qed
\end{corollary}

A leaf $\ell$ of a consistent separation tree on~$\vS$ is an \emph{$\F$-tangle leaf\/}, or {\em tangle leaf\/} for short, if the closure $\lfloor\beta_{\ell}\rfloor$ of $\beta_{\ell}$ is an $\F$-tangle of~$S$. A~non-leaf node~$v$ is a {\em tangle node\/} if there is a tangle leaf~$\ell\ge_r w$ for every successor~$w$ of~$v$ in~$>_r$.

A~leaf $\ell$ is \emph{forbidden} (by~$\F$) if $\beta_{\ell}$ contains an element of~$\F$ as a subset. Note that tangle leaves are never forbidden. 

\begin{definition}\label{def:TST}
An \emph{$\F$-tangle structure tree of~$\vS$} is a consistent separation tree on~$\vS$ in which every leaf is either a tangle leaf or forbidden, and for every non-leaf node~$v$ the set $\beta_v$ has no subset in $\F$.%
   \COMMENT{}
\end{definition}

Every $\F$-tangle structure tree of~$\vS$ displays all the $\F$-tangles of~$S$, not just some of them:

\begin{theorem}\label{thm:C}
    Let $\vS$ be a separation system, let $\F$ be any set, and let $T$ be any $\F$-tangle structure tree of~$\vS$. 
    \begin{enumerate}\itemsep2pt\vskip2pt
        \item\label{item:C1} For every $\F$-tangle $\tau$ of~$S$ there is a unique leaf $\ell$ of $T$ such that $\lfloor \beta_{\ell} \rfloor = \tau$. 
        \item\label{item:C2} If all the leaves of $T$ are forbidden, then $S$ has no $\F$-tangle.
    \end{enumerate}
\end{theorem}

\begin{proof}
\labelcref{item:C1} By \cref{lem:orichoosesleaf} there is a unique leaf $\ell$ with $\beta_{\ell} \subseteq \tau$. This leaf is not forbidden, since $\beta_{\ell}\subseteq\tau$ has no subset in~$\F$. As $T$ is a tangle structure tree, this means that $\ell$ is a tangle leaf: that $\lfloor \beta_{\ell} \rfloor$ is a tangle of~$S$. As $\lfloor \beta_{\ell} \rfloor \subseteq \tau$, this tangle can only be $\tau$. 

\labelcref{item:C2} is immediate from \labelcref{item:C1}.
\end{proof}

In the next section we shall determine the sets $\F$ for which $S$ admits a tangle structure tree whose tree-order is compatible with a given order function on~$S$, in that lower-order separations label edges further down on the tree. As we shall see in \cref{thm:Sk} and \cref{cor:tanglesinS}, such trees display not only all the $\F$-tangles of the entire set $S$, as all tangle structure trees do by \cref{thm:C}, but all the $\F$-tangles in~$\vS$: the $\F$-tangles of the sets $S_k = \{s \in S: |s| < k\}$ with $k \in \mathbb{N}$.

\medskip\section{Existence of tangle structure trees}\label{sec:TSTexistence}

\noindent
   Let $\vS$ be an ordered separation system, and let $\F$ be any~set. We shall need two conditions on~$\F$ to ensure that $\vS$ has a tangle structure tree.

The first is that $\{\sv\}\in\F$ for every $\vs \in \vS$ that is trivial in~$\vS$.  This holds for all sets~$\F$ of interest, and if it does we call ~$\F$ \emph{standard for~$\vS$}. Intuitively, co-trivial separations~$\sv$ point to places in our graph or structure that are tiny~-- too small to house any tangle. More formally, we already saw in \cref{sec:defs} that no consistent orientation of~$S$ can contain co-trivial separations. Assuming that $\F$ is standard for~$\vS$ therefore places no restrictions on the $\F$-tangles of~$S$.

The second condition has more substance. If $T$ is any separation tree on~$\vS$%
   \COMMENT{}
   then, by \cref{lem:orichoosesleaf}, every orientation~$\tau$ of~$S$ contains~$\beta_\ell$ for a unique leaf~$\ell$. If~$\tau$ is consistent and $\ucl(\beta_\ell)$ orients all of~$S$, then $\ucl(\beta_\ell) = \tau$ by \cref{lem:closure}.%
   \COMMENT{}
   If $\tau$ is not a tangle, then $\tau = \ucl(\beta_\ell)$ has a subset~$\sigma$ in~$\F$ that witnesses this.
If $T$ is in fact a tangle structure tree, we know that our consistent orientation~$\tau$ of~$S$ is a tangle (without having to assume that $\ucl(\beta_\ell)$ orients all of~$S$)%
   \COMMENT{}
   unless it contains such a set $\sigma\in\F$ not only in~$\ucl(\beta_\ell)$ but even in~$\beta_\ell$: among the edge labels of~$T$.

In order for this to be possible, we therefore need to make some `richness' assumption about~$\F$: an assumption which, in our example, ensures that $\F$ has enough elements to contain a subset also of~$\beta_\ell$ as soon as it contains a subset of~$\ucl(\beta_\ell)$.

We shall identify an essentially weakest-possible such richness condition below, in \cref{def:rich}. This will need some more preparation. For motivation, readers are invited to peek at \cref{minimization} and the definition preceding it now. That definition clearly implies that $\F$ contains a subset of~$\beta_\ell$ whenever it contains a subset of~$\ucl(\beta_\ell)$. The notion of richness from \cref{def:rich} will still imply this, but is weaker.

A~separation tree $T$ on an ordered separation system~$\vS$ is \emph{ordered} if $|s_v|\le |s_w|$ whenever $v$ and $w$ are non-leaves of $T$ with $v\le w$.%
   \COMMENT{}
   It is \emph{thoroughly ordered\/} (in~$\vS$) if, for every non-leaf node~$v$, the separation $s_v$ is not oriented by $\lfloor \beta_v \rfloor$%
   \COMMENT{}%
   \footnote{If $T$ is consistent, then this holds by \cref{cor:B}.}
   and has minimum order among the separations in~$S$ not oriented by~$\ucl(\beta_v)$.
 
\begin{lemma}\label{lem:A}
    Every thoroughly ordered separation tree is ordered.
\end{lemma}

\begin{proof} Let~$T$ be a thoroughly ordered separation tree on~$\vS$. If it is not ordered, it has nodes $v, w$ with $v < w$ such that $|s_v| > |s_w|$. Since $s_w$ is not oriented by $\lfloor \beta_v \rfloor \subseteq \lfloor \beta_w \rfloor$,%
   \COMMENT{}
   this contradicts the requirement that $s_v$ have minimum order among the separations not oriented by $\lfloor \beta_v \rfloor$. 
\end{proof}

We say that $\vs \in \vS$ is \emph{weakly eclipsed} by $\vr \in \vS$ if $\vr < \vs$ and $|r|\le |s|$, and \emph{eclipsed} by~$\vr$ if $\vr<\vs$ and $|r| < |s|$. Given any set $\tau \subseteq \vS$, a subset $\sigma \subseteq \tau$ is \emph{efficient} (in~$\tau$) if no element of $\sigma$ is eclipsed by any other element of~$\tau$. It is \emph{strongly efficient} if no element of $\sigma$ is weakly eclipsed by any other element of~$\tau$. Note that if the order function on~$S$ is injective and~$\tau$ is a partial orientation of~$S$, then every efficient subset of~$\tau$ is strongly efficient in~$\tau$.%
   \COMMENT{}%
   \footnote{We need the assumption on~$\tau$, since $\vs$ eclipses~$\sv$ weakly if $\vs < \sv$. So we do not want $\vs,\sv\in\tau$.\looseness=-1}

\begin{lemma}\label{lem:thorough}
Let $v$ be a node of a thoroughly ordered separation tree~$T\!$ on~$\vS$. Then
\begin{enumerate}\itemsep2pt\vskip2pt
    \item\label{efficient_and_ordering_1} every strongly efficient subset of $\lfloor \beta_v \rfloor$ is contained in~$\beta_v$;
    \item\label{efficient_and_ordering_2} if $\beta_v$ is consistent, it is efficient in any partial orientation $\tau$ of~$S$ that includes~$\ucl(\beta_v)$.%
   \COMMENT{}
\end{enumerate}
\end{lemma}

\begin{proof}
\labelcref{efficient_and_ordering_1}  Let $\sigma$ be any strongly efficient subset of $\lfloor \beta_v \rfloor$, and let $\vs\in\sigma$ be given. 
Suppose that $\vs\notin\beta_v$. Then $\vs \in \lfloor \beta_v \rfloor \setminus \beta_v$ is required by $\beta_v$, so there exists $\vr\in\beta_v$ with $\vr<\vs$ (and $r\neq s$). As $\vr\in\beta_v$, we have $r=s_u$ for some $u<v$. Choose such an $\vr$ with $u$ minimal in~$<_r$.

Our aim is to show that $|r|\le |s|$: then $\vr \in \lfloor \beta_v\rfloor$ eclipses $\vs \in \sigma$ weakly, contradicting the strong efficiency of $\sigma$ as a subset of $\lfloor \beta_v \rfloor$. 
Since $T$ is thoroughly ordered, we shall have $|r|\le |s|$ as desired if $\lfloor \beta_u \rfloor$ does not orient~$s$: this would make~$s$ a candidate for~$s_u$, so $|r|>|s|$ would contradict the fact that $r=s_u$.%
   \COMMENT{}

So let us show that neither $\vs$ nor~$\sv$ lies in~$\lfloor \beta_u \rfloor$. For $\vs \notin \lfloor \beta_u \rfloor$, recall first that $\vs\notin\beta_v\supseteq\beta_u$. Hence if $\vs\in\ucl(\beta_u)$, there exists an $\vrdash\in\beta_u$ such that $\vrdash < \vs$. This satisfies $r'=s_{u'}$ for some~$u'<_r u$, so $\vrdash < \vs$ contradicts our original choice of $\vr$ given~$\vs$. Hence $\vs\notin\ucl(\beta_u)$ as desired. Suppose now that $\sv \in \lfloor \beta_u \rfloor$. As $\vs > \vr$, any $\vrdash \leq \sv$ in $\beta_u$ satisfies $\vrdash\leq \sv < \rv$. So $r$ was already oriented by $\lfloor \beta_u \rfloor$, contradicting the fact that $r=s_u$.%
   \COMMENT{}
   Hence neither $\vs$ nor~$\sv$ lies in~$\lfloor \beta_u \rfloor$, as desired.

\labelcref{efficient_and_ordering_2} Suppose not; then some $\vr\in\tau$ eclipses some $\vs\in \beta_v$. 
Let $u< v$ be such that $s = s_u$. As $|r| < |s|$ and $T$ is thoroughly ordered, we know that $\lfloor \beta_u \rfloor$ orients $r$.

Since $\tau$ is a partial orientation of~$S$ containing $\vr$, it cannot also contain $\rv$. We thus cannot have $\rv \in \lfloor \beta_u \rfloor$, since $\lfloor \beta_u \rfloor \subseteq \lfloor \beta_v \rfloor \subseteq \tau$ by \cref{lem:closure}.%
   \COMMENT{}
   So $\vr \in \lfloor \beta_u \rfloor$. As $\vr < \vs$,%
   \COMMENT{}
   this implies $\vs\in\ucl(\beta_u)$ unless $\vr = \sv\in\beta_u$.%
   \COMMENT{}
   Both these contradict the fact that $s = s_u$. 
\end{proof}

A~typical application of \cref{lem:thorough}\,\ref{efficient_and_ordering_2} is that $\beta_v$ is efficient in any tangle of~$S$ that contains it. Such a tangle will be consistent and hence include~$\ucl(\beta_v)$ by \cref{lem:closure}.%
   \COMMENT{}

A typical application of \cref{lem:thorough}\,\ref{efficient_and_ordering_1} will be that witnesses $\sigma\in\F$ to the fact that some $\ucl(\beta_v)$ fails to extend to a tangle can be found not only in~$\ucl(\beta_v)$ but in~$\beta_v$ itself. Then our structure trees will display such witnesses $\sigma\in\F$ in the label sets~$\beta_\ell$ of their forbidden leaves~$\ell$. For this to work with the help of \cref{lem:thorough}, we need $\F$ to contain witnesses that are strongly efficient in their~$\ucl(\beta_v)$.

This motivates our formal richness condition on~$\F$:

\begin{definition}\label{def:rich}
A set $\F$ is \emph{rich for~$\vS$} if every consistent orientation of~$S$ that has a subset in~$\F$ also has a strongly efficient\footnote{in this orientation of~$S$} subset in~$\F$.%
   \COMMENT{}
\end{definition}

The assumption that $\F$ is rich will be central to our results. It is needed to ensure that our $\F$-tangle structure trees exist, and we shall see after \cref{lemma-twoooo} that no weaker condition on~$\F$ will ensure the same.

Although it may look a bit technical at first glance, this `richness' requirement on~$\F$ is quite natural, given the role of these~$\F$ in tangle theory.
For example, given a consistent orientation~$\tau$ of~$S$ that has a subset $\sigma\in\F$, we can often obtain a strongly efficient subset~$\sigma'\in\F$ of~$\tau$ simply by replacing every element~$\vs$ of~$\sigma$ by some $\vsdash\le\vs$ that is minimal in~$\tau$.%
   \COMMENT{}
   One still has to check then that $\sigma'$ is indeed in~$\F$. But as the idea behind those forbidden triples in~$\F$ is that they identify areas in our graph or other structure that are `too small to be home to a tangle', it is not unnatural for this particular~$\sigma'$ to be in~$\F$ if~$\sigma$ was.

Let us cast this example in the form of a lemma. Let us call a set~$\F$ {\em closed under minimization\/} in a subset~$\tau$ of~$\vS$ if it contains every set~$\sigma'\subseteq\tau$ obtained from some $\sigma\subseteq\tau$ in~$\F$ by replacing every $\vs\in\sigma$ with some $\vsdash\le\vs$.

\begin{lemma}\label{minimization}
If $\F$ is closed under minimization in every consistent orientation of~$S$,%
   \COMMENT{}
   then $\F$ is rich for~$\vS$.
\end{lemma}%
   \COMMENT{}

\begin{proof}
Let $\tau$ be any consistent orientation of~$S$ that has a subset~$\sigma$ in~$\F$. We have to find a set $\sigma'\subseteq\tau$ in~$\F$ that is strongly efficient in~$\tau$.

Let $\sigma'$ be obtained from~$\sigma$ by replacing every $\vs\in\sigma$ by some $\vsdash\le\vs$ that is minimal in~$\tau$.%
   \COMMENT{}
   The set $\sigma'$ is strongly efficient in~$\tau$, since any $\vr < \vsdash\in\sigma'$ in~$\tau$ contradicts the minimal choice of $\vsdash\le\vs$. Since $\F$ is closed under minimization in~$\tau$, we have $\sigma'\in\F$ as required.
\end{proof}

When we apply \cref{minimization} later in \cref{sec:apps}, we shall in fact prove that the sets~$\F$ needed there are closed under minimization in all of~$\vS$. Note that this is stronger than being closed under minimization in every consistent orientation of~$S$.

\begin{lemma}\label{coreidea}
Let $\F$ be rich and standard for~$\vS$. Let $T$ be a thoroughly ordered separation tree on~$\vS$, and let $v \in V(T)$. Assume that $\beta_v$ is consistent and avoids~$\F$, and that $\lfloor \beta_v \rfloor$ orients all of~$S$. Then $\lfloor\beta_v\rfloor$ is a tangle of~$S$.
\end{lemma}

\begin{proof}
As $\beta_v$ avoids $\F$, which is standard, $\beta_v$ has no element that is co-trivial in~$\vS$. As $\beta_v$ is consistent, $\lfloor\beta_v\rfloor$~is consistent by \cref{lem:closure-consistent}. 
Since $\lfloor \beta_v \rfloor$ orients all of~$S$, it can thus only fail to be a tangle of~$S$ if it has a subset~$\sigma$ in~$\F$. As $\F$ is rich for~$\vS$, we can choose $\sigma$ to be strongly efficient as a subset of $\lfloor \beta_v \rfloor$.%
   \footnote{If $\lfloor \beta_v \rfloor$ contains both orientations of some $s \in S$, delete one of them to obtain a consistent orientation of~$S$ as required in the definition of `rich' for $\F$.}
   Then even $\sigma \subseteq \beta_v$ by
\cref{lem:thorough}\,\labelcref{efficient_and_ordering_1}, which contradicts our assumption that $\beta_v$ avoids $\F$. Hence $\lfloor\beta_v\rfloor$ is a tangle of~$S$.
\end{proof}

\begin{theorem}\label{TreeConstruction}
Let $\vS$ be an ordered separation system, and let $\F$ be a set that is rich and standard for~$\vS$. 
Then there exists a thoroughly ordered $\F$-tangle structure tree of~$\vS$.
\end{theorem}

\begin{proof}
   Start with the one-node tree $T_0=\{r\}$ and let $\beta^0 = \emptyset$. We iteratively build trees $T_0\subsetneq T_1\subsetneq\dots$ with maps $\beta^i$ so that each $(T_i,r,\beta^i)$ is consistent and thoroughly ordered, and none of their sets~$\beta^i_v$ with $v$ a non-leaf node will have a subset in~$\F$. Note that $T_0$ has all these properties. The last of those trees will be our desired tangle structure tree.

If for some $n$ the tree $(T_n,r,\beta^n)$ is already a tangle structure tree, we are done. Otherwise $T_n$ has a leaf~$v$ which is neither a tangle leaf nor forbidden. By \cref{coreidea}, $\lfloor\beta^n_v\rfloor$ does not orient all of~$S$; let $s\in S$ be a%
   \COMMENT{}
   separation of minimum order not oriented by~$\lfloor\beta^n_v\rfloor$.%
   \COMMENT{}

Form $T_{n+1}$ by adding two children $v_1,v_2$ at~$v$. Let $\beta^{n+1}$ agree with~$\beta^n$ on the edges of~$T_n$, and pick orientations $\vs =: \beta^{n+1}(vv_1)$ and $\sv =: \beta^{n+1}(vv_2)$ of~$s$; then $s = s_v$ in~$T_{n+1}$. By construction, $(T_{n+1},r,\beta^{n+1})$ is a thoroughly ordered separation tree. It is consistent, because $T_n$ was and neither orientation of~$s$ lies in $\lfloor\beta^{n+1}_v\rfloor$. Its only non-leaf node~$v$ that was not already a non-leaf node in~$T_n$ is~$v$. Since $v$ was not forbidden as a leaf of~$T_n$, the set~$\beta^{n+1}_v = \beta^n_v$ has no subset in~$\F$.

This process strictly increases $|V(T_n)|$ at each step, but by \cref{lem:smoll} there is an upper bound on the size of any consistent separation tree in terms of~$|S|$. Hence the process terminates after finitely many steps with a thoroughly ordered tangle structure tree.
\end{proof}

If our order function on~$S$ is injective, \cref{TreeConstruction} has a converse, which shows that our requirement of richness for~$\F$ is weakest possible to ensure the existence of a tangle structure tree. This is established by our next lemma:

\begin{lemma}\label{lemma-twoooo}
Let $\vS$ be a separation system with an injective order function on~$S$, and let~$\F$ be any set. 
If there exists a thoroughly ordered%
   \COMMENT{}
   $\F$-tangle structure tree~$T$ of~$\vS$, then $\F$ is rich for~$\vS$.
\end{lemma}

\begin{proof}
Let $\tau$ be any consistent orientation of~$\vS$ that has a subset in $\F$. We shall find an efficient subset $\sigma \in \F$ of~$\tau$, which will even be strongly efficient since our order function is injective.%
   \COMMENT{}
   By \cref{lem:orichoosesleaf}, $T$~has a leaf $\ell$ with $\beta_\ell\subseteq\tau$. In particular, $\beta_{\ell}$ is consistent, so $\lfloor \beta_{\ell} \rfloor \subseteq \tau$ by \cref{lem:closure}.

Since $T$ is a tangle structure tree, $\ell$ is either a tangle leaf or forbidden. If $\ell$ is a tangle leaf then $\lfloor\beta_\ell\rfloor$ is a tangle of all of~$S$. But then $\lfloor \beta_{\ell} \rfloor = \tau$, which contradicts the fact that $\tau$ has a subset in~$\F$. 

Thus $\ell$ is forbidden, so $\beta_\ell$ has a subset $\sigma$ in $\F$. By 
\cref{lem:thorough}\,\labelcref{efficient_and_ordering_2} this~$\sigma$ is efficient in~$\tau$, as desired.
\end{proof}

\noindent
   Note that if our order function on~$S$ is not injective, the proof of \cref{lemma-twoooo} still goes through as stated except for one aspect: the efficient subset $\sigma\in\F$ of~$\tau$ it finds may not be strongly efficient in~$\tau$ (as our definition of `rich' requires).

\medbreak

If our order function on~$S$ is injective, the sets~$\F$ for which $\vS$ admits a thoroughly ordered tangle structure tree are thus precisely the rich ones:

\begin{theorem}\label{thm:generalduality}
Let $\vS$ be a separation system with an injective order function on~$S$, and let~$\F$ be a set that is standard for~$\vS$. 
There exists a thoroughly ordered $\F$-tangle structure tree of~$\vS$ if and only if $\F$ is rich for~$\vS$. If one exists, it is unique.
\end{theorem}

\begin{proof}
Existence is immediate from \cref{TreeConstruction,lemma-twoooo}. Uniqueness follows from the definition of `thoroughly ordered' and our assumption that the order function on~$S$ is injective.
\end{proof}

Let us illustrate the above by an example. Consider a separation system whose elements are the non-empty subsets of a 4-element set~$V\!$, with $\le$ defined as~$\subseteq$, and $^*$ as complementa\-tion in~$V\!$. Let $r,s$ be two crossing separations; their orientations are 2-sets. The 1-element subsets of~$V\!$ then are precisely those of the form $\vr\cap\vs$ for suitable orientations $\vr,\vs$ of $r$ and~$s$. We call the four unoriented separations which each partition~$V\!$ into a singleton and a 3-set the {\em corners\/} of $r$ and~$s$. Let $S$ be the set of $r$, $s$, and their four corners.

The set~$S$ then has four {\em principal\/} consistent orientations, those that orient one corner~$t$ away from $r$ and~$s$ and all the other elements of~$S$ towards~$t$. And it has the {\em non-principal\/} consistent orientations, which orient $r$ and $s$ arbitrarily and all four corners towards $r$ and~$s$. Choose an injective order function on $S$ such that $|r| < |s| < |t|$ for all corners~$t$.

Let us first consider as~$\F$ the 4-element%
   \COMMENT{}
   set $\P = \big\{\{\vr,\vs,\!\tv\} \mid \vt = \vr\cap\vs\big\}$.%
   \footnote{We shall study these $\P$-tangles more closely in \cref{sec:apps}.}
   The unique thoroughly ordered separation tree~$T$ on~$\vS$ then starts with orienting $r$ and~$s$, as a tree~$T'$ with four leaves, one for every orientation of~$\{r,s\}$. For each of these leaves~$v$, the next separation~$s_v$ to be oriented is one of the corners~$t$. If $\vt = \vr\cap\vs$, say,%
   \COMMENT{}
   we have $\vt = \beta(v\ell)$ for a tangle leaf~$\ell$ for the principal tangle~$\ucl(\beta_\ell) = \ucl(\vt)$. And $\!\tv = \beta(v\ell')$ for a forbidden leaf~$\ell'$, with $\beta_{\ell'} = \{\vr,\vs,\!\tv\}\in\P$. As this happens at each of the four leaves~$v$ of~$T'\!$, our $T$~is a tangle structure tree of~$\vS$ for $\F = \P$.\looseness=-1

Now let~$\F$ be the set~$\RR$ of all triples in~$\vS$ whose complements partition~$V\!$.%
   \footnote{We shall study these $\RR$-tangles more closely in~\cite{TSTs2}.}
   These are the triples of the form either $\{\vr,\!\tv,\tvdash\}$ with $\vt = \vr\cap\vs$ and $\vtdash = \vr\cap\sv$, or $\{\vs,\!\tv,\tvdash\}$ with $\vt = \vs\cap\vr$ and $\vtdash = \vs\cap\rv$. Our thoroughly ordered separation tree~$T$ then starts with the same~$T'$ as before. Every leaf~$v$ of~$T'$ once more sends an edge $v\ell$ to a tangle leaf~$\ell$ with $\beta(v\ell) = \!\vt$ and $\RR$-tangle $\ucl(\beta_\ell) = \ucl(\vt)$. But with $\beta_v = \{\vr,\vs\}$ as earlier we now have $\beta_{\ell'} = \{\vr,\vs,\!\tv\}\notin\RR$ for the other child~$\ell'$ of~$v$. But note that $\tvdash\supset\vs$ for the corner $\vtdash = \vr\cap\sv$. Hence $\ucl(\beta_{\ell'})$ is an orientation of~$S$ with a subset in~$\RR$: the set $\{\vr,\!\tv,\tvdash\}$. Thus, $\ell'$~is neither a tangle leaf nor forbidden, and $T$ is not an $\RR$-tangle structure tree.

As predicted by \cref{thm:generalduality}, this difference between $\P$ and~$\RR$ is reflected by the fact that $\P$ is rich for~$\vS$ but $\RR$ is not: the consistent orientation~$\ucl(\beta_{\ell'})$ of~$S$ has a subset in~$\RR$, but it has no efficient subset in~$\RR$. In particular, its subset $\{\vr,\!\tv,\tvdash\}\in\RR$, which prevented~$T$ from being an $\RR$-tangle structure tree, is not efficient, because $\vs\subset\tvdash$ eclipses~$\tvdash$.

\medbreak

In \cref{sec:efficient} we shall prove that our tangle structure trees can be improved further, without any additional requirements on~$\F$, by contracting `inessential' edges: edges whose label~$\vs$ is neither needed in any sets $\sigma\in\F$ witnessing that a leaf is forbidden, nor needed to determine any tangles.%
   \COMMENT{}

Unfortunately, this contraction process can cause our structure trees to lose their property of being thoroughly ordered. We shall therefore extract from this property its essence, a~slightly weaker property to be called `efficiency', which is still strong enough to make the branches~$rT\ell$ and their label sets~$\beta_\ell$ as efficient for displaying the tangles~$\ucl(\beta_\ell)$ as they are in thoroughly ordered structure trees, but weak enough to survive the contraction process we envisage for \cref{sec:efficient}.

\medbreak

Given an ordered separation system~$\vS$, we call a separation tree $(T,r,\beta)$ \emph{efficient\/}%
   \COMMENT{}
   if for every leaf~$\ell$ the set~$\beta_\ell$%
   \COMMENT{}
   is efficient in~$\lfloor \beta_\ell \rfloor$.%
   \COMMENT{}

\begin{lemma}\label{lem:E}
Let $\vS$ be an ordered separation system. Then every thoroughly ordered separation tree $(T,r_0,\beta)$ on~$\vS$ is efficient.
\end{lemma}

\begin{proof}
   Consider any leaf~$\ell$ of~$T$.  If $\beta_\ell$ is not efficient in~$\ucl(\beta_\ell)$ as claimed, then some $\vs\in\beta_\ell$ is eclipsed by some $\vr\in\ucl(\beta_\ell)$. Let $v,w<\ell$ be such that $\vs=\vsv$ and $\vr\ge\vsw\in\beta_\ell$.%
   \footnote{Recall that our arrow notation for oriented separations is never fixed. We are thus free to use forward arrows to denote the orientations of $s_v$ and~$s_w$  that lie in $\beta_\ell$.}

As $\vr$ eclipses~$\vs$, we have $\vr < \vs$ and $|r| < |s|$;%
   \COMMENT{}
   in particular, $r\ne s=s_v$. The fact that $T$ is thoroughly ordered thus implies that $r$ is already oriented by~$\ucl(\beta_v)$. We cannot have $\vr\in\ucl(\beta_v)$, since this would place $\vsv = \vs > \vr$ or~$\svv$ in~$\ucl(\beta_v)$ too,%
  \footnote{By the definition of~$\ucl(\ )$, the assumption of $\vs > \vr\in\ucl(\beta_v)$ implies $\vs\in\ucl(\beta_v)$ only if $\sv\notin\beta_v$.}
   which would contradict the fact that $T$ is thoroughly ordered.%
   \COMMENT{}
   Thus, $\rv\in\ucl(\beta_v)$; let $u < v$ be such that $\rv\ge\vsu\in\beta_v$.

If $u=w$, then%
   \COMMENT{}
   $\vsu = \vsw\le\vr < \vs$ and hence $\vsv = \vs\in\ucl(\beta_v)$ or $\svv\in\beta_v$,$^{\thefootnote}$ both of which contradict the fact that $T$ is thoroughly ordered.%
   \COMMENT{}

If $u < w$, then $\svw\ge\rv\ge\vsu\in\beta_w$.%
   \COMMENT{}
   Then $\svw\in\ucl(\beta_w)$ unless $\vsw\in\beta_w$,$^{\thefootnote}$ which both contradict the fact that $T$ is thoroughly ordered.%
   \COMMENT{}

If $w < u$, finally, then $\svu\ge\vr\ge\vsw\in\beta_u$,%
   \COMMENT{}
   so $\svu\in\ucl(\beta_u)$ unless $\vsu\in\beta_u$,$^{\thefootnote}$ which both contradict the fact that $T$ is thoroughly ordered.
   \end{proof}%
   \COMMENT{}

\section{Irreducible and efficient tangle structure trees}\label{sec:efficient}

\noindent
   Let $(T,r,\beta)$ be a consistent%
   \COMMENT{}
   separation tree on a separation system~$\vS$. Let $v$ be a node of~$T$ with a child~$w$. 

Let us call the edge $vw$ of~$T$ \emph{necessary\/} for a tangle leaf $\ell\ge w$ of~$T$ if $\beta(vw)$ is a $\le$-minimal element of the tangle~$\ucl(\beta_\ell)$, or equivalently, of~$\beta_\ell$. We call~$vw$ \emph{necessary\/} for a forbidden leaf~$\ell\ge w$ if every subset of~$\beta_\ell$ in~$\F$ contains~$\beta(vw)$.

Let $(T,r,\beta)_{w\to v} := (T_{w\to v},\, r_{w\to v},\, \beta_{w\to v})$,
where $T_{w \to v}$ is obtained from $T$ by contracting the edge $vw$ of~$T$ and deleting any other child~$w'$ of~$v$ together with the subtree spanned by the nodes $u\ge_r w'$.
We continue to use `$w$' for the new node constructed from the edge $vw$, and thus think of $V(T_{w \to v})$ as a subset of $V(T)$.  Similarly, we think of $E(T_{w \to v})$ as a subset of $E(T)$ and continue to use `$\beta$' to denote its labelling $\beta_{w \to v}$.  If $v=r$ we set $r_{w\to v}:=w$; otherwise we keep $r_{w\to v}:=r$.\looseness=-1

Note that $T_{w \to v}$ inherits the sets $E_u$ from $T$ for its nodes~$u$. In particular, leaves of~$T_{w\to v}$ were also leaves of~$T$, and the separations $s_u$ for non-leaf nodes~$u$ of~$T_{w\to v}$ are still what they were for~$T$. Hence $T_{w \to v}$ is still ordered if $T$ was. However, $T_{w \to v}$ may no longer be thoroughly ordered in~$\vS$ even if $T$~was.

The sets $\beta_u$ remain unchanged for all $u \ngeq w$, while for all $u\ge w$ the sets $\beta_u$ in~$T_{w \to v}$ arise from $\beta_u$ in~$T$ by deleting~$\beta(vw)$ from it. (Recall that $\beta$ was injective on~$E(rTu)$, by \cref{def:septree}.) In particular, $T_{w \to v}$ is consistent if $T$~was. For tangle leaves~$\ell$ of~$T_{w\to v}$, the tangle~$\ucl(\beta_\ell)$ of~$S$ is the same as it was when $\beta_\ell$ was taken in~$T$; its former element~$\beta(vw)$ then was a non-minimal element of this tangle.

\begin{lemma}\label{lem:necessity}
    Let $\vS$ be a separation system, and let~$\F$ be any set. Let $(T, r_0, \beta)$ be a consistent separation tree%
   \COMMENT{}
   on~$\vS$. Let $v$ be a node of~$T$ with a successor~$w$, let $\vs:=\beta(vw)$, and let $\ell\ge w$ be a leaf. Then $vw$ is necessary for~$\ell$ if and only if $\beta_{\ell} \setminus \{\vs\}$ has no subset in~$\F$ and $\lfloor \beta_{\ell} \setminus \{\vs\} \rfloor$ is not a tangle of~$S$.\looseness=-1
\end{lemma}

\begin{proof}
Suppose first that $vw$ is necessary for~$\ell$, and that $\ell$ is a tangle leaf. Then $\ucl(\beta_\ell)$ is a tangle of~$S$, so $\beta_{\ell} \setminus \{\vs\}$ has no subset in~$\F$.

For a proof that $\lfloor \beta_{\ell} \setminus \{\vs\} \rfloor$ is not a tangle of~$S$, we show that it contains neither $\vs$ nor~$\sv$. To see $\vs\notin\lfloor \beta_{\ell} \setminus \{\vs\} \rfloor$, note that any $\vr\in\beta_\ell\setminus\{\vs\}$ with $\vr\le\vs$%
   \COMMENT{}
   would satisfy~$\vr<\vs$, since also $\vs\in\beta_\ell$ and hence $r\ne s$ by \cref{def:septree}. Then $\vs$ would not be minimal in~$\beta_\ell$, which it is since $vw$ is necessary for the tangle leaf~$\ell$. And $\sv\notin\lfloor \beta_{\ell} \setminus \{\vs\} \rfloor$, since $\lfloor \beta_\ell \rfloor$ is a tangle of~$S$ containing~$\vs$. Thus, $\lfloor \beta_{\ell} \setminus \{\vs\} \rfloor$ does not orient~$s$, so it is not a tangle of~$S$. 

Suppose, second, that $vw$ is necessary for~$\ell$, and that $\ell$ is a forbidden leaf. Then $\beta_{\ell} \setminus \{\vs\}$ has no subset in~$\F$. Suppose $\lfloor \beta_{\ell} \setminus \{\vs\} \rfloor$ is a tangle of~$S$. This tangle cannot contain~$\vs$: it would then contain the entire set~$\beta_\ell$, but this has a subset in~$\F$ since $\ell$ is a forbidden leaf. But neither can $\lfloor \beta_{\ell} \setminus \{\vs\}\rfloor$ contain $\sv$: then $\beta_{\ell} \setminus \{\vs\}$ would contain some $\vr \leq \sv$, which would contradict the consistency of~$\beta_{\ell}$, since $r \neq s$ by \cref{def:septree}%
   \COMMENT{}
   and $\vs \in \beta_{\ell}$ would thus point away from $\vr \in \beta_{\ell}$. 

Suppose, third, that $vw$ is not necessary for~$\ell$ and that $\ell$ is a tangle leaf. Then $\ucl(\beta_\ell)$ is a tangle and $\vs$ is not a minimal element of~$\beta_\ell$. We complete our proof in this case by showing that $\ucl(\beta_\ell\setminus\{\vs\}) = \ucl(\beta_\ell)$.%
   \COMMENT{}
   If $\vs\notin\beta_\ell$ this holds trivially. If $\vs$ lies in~$\beta_\ell$ but is not minimal in it, then there exists some $\vr\in\beta_\ell$ with $\vr < \vs$. This implies $\vs\in\ucl(\beta_\ell\setminus\{\vs\})$, and hence
 $$\ucl(\beta_\ell) = \ucl(\beta_\ell\setminus\{\vs\})\cup \ucl(\{\vs\})
   \subseteq \ucl(\ucl(\beta_\ell\setminus\{\vs\})) = \ucl(\beta_\ell\setminus\{\vs\})$$
 as desired, unless $\sv\in\ucl(\beta_\ell\setminus\{\vs\})$. But in that case the tangle~$\ucl(\beta_\ell)$ contains both $\vs$ and~$\sv$, which it cannot.

Suppose finally that $vw$ is not necessary for~$\ell$ and that $\ell$ is a forbidden leaf. Then $\beta_{\ell} \setminus \{\vs\}$ has a subset in~$\F$, as desired. 
\end{proof}

\begin{lemma}\label{lem:edge_contraction}
Let $\vS$ be a separation system, and let~$\F$ be any set. Let $(T,r,\beta)$ be a tangle structure tree of~$\vS$. Let $v$ be a node of~$T$ with a successor~$w$. Then $(T,r,\beta)_{w\to v}$ is a tangle structure tree of~$\vS$ if and only if $vw$ is not necessary for any leaf $\ell\ge w$ of~$T$.
\end{lemma}

\begin{proof}
Suppose first that $vw$ is necessary for some leaf $\ell\ge w$ of~$T$.  This $\ell$ is a leaf also of~$T_{w \to v}$. By \cref{lem:necessity}, however, $\ell$~is neither a tangle leaf nor a forbidden leaf of~$T_{w \to v}$.%
   \COMMENT{}
   Hence $T_{w\to v}$ is not a tangle structure tree of~$\vS$.

Conversely, suppose that $vw$ is not necessary for any leaf $\ell \geq w$ in~$T$. These~$\ell$ are precisely (all) the leaves of~$T_{w \to v}$. Since $\beta_\ell$ taken in~$T_{w \to v}$ is obtained from $\beta_\ell$ taken in~$T$ by deleting~$\beta(vw)$,%
   \COMMENT{}
   these~$\ell$ are either forbidden or tangle leaves of~$T_{w \to v}$ by \cref{lem:necessity}. Thus, $T_{w \to v}$ is a tangle structure tree of~$\vS$.%
   \COMMENT{}
   \end{proof}

Let us note a surprising consequence of \cref{lem:edge_contraction} and its proof. Consider a node~$v$ of~$T$ with children~$w,w'$, and assume that $vw$ is not necessary for any leaf $\ell\ge w$ of~$T$. The truth of this assertion is determined entirely by the values under~$\beta$ of the subtree of~$T$ that consists of the path~$rTv$, the subtree spanned by all nodes~$u\ge w$, and the edge~$vw$. In particular, it is not affected by the values under~$\beta$ of the edges of the subtree of~$T$ spanned by the nodes~$u\ge w'$: the subtree we delete when we form $T_{w \to v}$ from~$T$.%
   \COMMENT{}

By \cref{lem:edge_contraction}, this $T_{w \to v}$~is a tangle structure tree. Every tangle of~$S$ thus equals~$\ucl(\beta_\ell)$ for some leaf~$\ell$ of~$T_{w \to v}$, which we recall is the same tangle of~$S$ that was associated with~$\ell$ in~$T$. (In particular, $\ucl(\beta_\ell)$ is the same in~$T_{w \to v}$ as in~$T$.) So what about the tangles of~$S$ whose associated leaf of~$T$ lay in the deleted subtree, tangles that include~$\beta_{w'}$ and in particular the inverse~$\beta(vw')$ of~$\beta(vw)$?

The surprising conclusion is that such tangles cannot exist: that if $vw$ is not necessary for any leaf~$\ell\ge w$, then $\beta_{w'}$ does not extend to a tangle of~$S$.%
   \footnote{This does not mean that $\beta(vw')$ cannot lie in a tangle of~$S$: such tangles just do not include~$\beta_v$.}
   This is true regardless of whether $\beta_w$ extends to a tangle of~$S$ or not.%
   \COMMENT{}
   But in particular, $v$~cannot have been a tangle leaf of~$T$; indeed, $T$~had no tangle node~$v'\ge w'$:

\begin{corollary}\label{cor:TangleNodesInContractions}
Tangle structure trees%
   \COMMENT{}
   $T\!$ and~$T_{w\to v}$ as in \cref{lem:edge_contraction} have the same%
   \COMMENT{}
   sets of tangle nodes.
\end{corollary}

It is not hard to construct examples of tangle structure trees that have unnecessary edges.%
   \COMMENT{}
   Indeed in many naturally occurring tangle structure trees most edges are unnecessary, and applying \cref{thm:irreducible} below reduces them substantially.%
   \COMMENT{}

\medbreak

Let us call a node~$v$%
   \COMMENT{}
   of~$T$ \emph{necessary in~$T$} if for every child $w$ of $v$ there exists a leaf $\ell \geq w$ of~$T$ such that $vw$ is necessary for $\ell$.

\begin{lemma}\label{lem:necessary_and_contraction}
Let $\vS$ be a separation system, and let~$\F$ be any set. Let $(T,r,\beta)$ be a tangle structure tree of~$\vS$, and let $v$ be a node of~$T$. Then $v$ has a child~$w$ such that $(T,r,\beta)_{w\to v}$ is a tangle structure tree of~$\vS$ if and only if $v$ is not necessary in~$T\!$.\looseness=-1
\end{lemma}

\begin{proof}
Suppose first that $v$ is necessary in~$T$. Let $w$ be any child of~$v$ in~$T$.%
   \COMMENT{}
   Since $v$ is necessary, there is a leaf $\ell\ge w$ such that $vw$ is necessary for~$\ell$. By \cref{lem:edge_contraction}, $T_{w\to v}$ is not a tangle structure tree of~$\vS$.

Conversely, suppose $v$ is not necessary in $T$. Then it has a child $w$ such that $vw$ is not necessary for any leaf $\ell \geq w$ in $T$. By \cref{lem:edge_contraction}, $T_{w \to v}$ is a tangle structure tree of~$\vS$.
\end{proof}

A tangle structure tree $(T,r,\beta)$ is called \emph{irreducible} if every node of~$T$ is necessary. The following theorem summarizes how tangle structure trees can be contracted and pruned to irreducible trees that inherit many of their properties. The statement of the theorem assumes our notational convention that $V(T_{w\to v})\subseteq V(T)$.

\begin{theorem}\label{thm:irreducible}
Let $\vS$ be a separation system, and let~$\F$ be any set.%
   \COMMENT{} 
For every $\F$-tangle structure tree~$T$ of~$\vS$ there exists an irreducible $\F$-tangle structure tree $T'$ of~$\vS$ obtained from $T$ by a sequence $T=T^1, \dots, T^n = T'$ with $T^{i+1} = T^i_{w \to v}$ for suitable nodes $v, w$ of~$T^i$.

The tree~$T'$ imposes the same partial ordering on its node as $T\!$ does.%
   \COMMENT{}
   In particular, $T'\!$~is ordered if $T\!$ is. Its leaves~$\ell$ are also leaves of~$T\!$, and their incident edges in~$T$ are still edges of~$T'\!$.

The tangle leaves~$\ell$ of~$T'\!$ are precisely those of~$T$.%
   \COMMENT{}
   Their associated tangles~$\ucl(\beta_\ell)$ are the same in~$T'\!$ as in~$T$. The tangle nodes of~$T'$ are precisely those of~$T$. A~tangle node~$v$ separates two tangle leaves in~$T'\!$ if and only if it separates them in~$T\!$.
\end{theorem}

\begin{proof}
If $T^i$ is not irreducible, pick a node $v \in T^i$ that is not necessary. By \cref{lem:necessary_and_contraction} it has a child $w$ for which $T^i_{w\to v}$ is a tangle structure tree. Let $T^{i+1} := T^i_{w\to v}$. For the second paragraph of the theorem, note that edges of the form~$v\ell$, with $\ell$ a leaf, are necessary in every tangle structure tree: by \cref{cor:B} if $\ell$ is a tangle leaf, and by \cref{def:TST} if $\ell$ is forbidden.%
   \COMMENT{}

For the third paragraph, recall our discussion following \cref{lem:edge_contraction} and iterate \cref{cor:TangleNodesInContractions}. Recall that the separation~$s_v$ associated with a tangle node~$v$ is the same in $T'\!$ as in~$T$, and that $s_v$ distinguishes two tangles if and only if $v$ separates their tangle leaves in any tangle structure tree of~$\vS$.
\end{proof}

The irreducible tangle structure tree~$T'$ obtained in the proof of \cref{thm:irreducible} is not unique. This is because we can check the nodes~$v$ of~$T$ for suppression in any order, and this order will affect their necessity in the tree~$T^i$ in which they are considered for suppression. But it is easy to make $T'$ unique. For example, it is not hard to show that if we consider the nodes~$v$ in any linear order that extends their tree order~$<_r$ from~$T$~-- for example, level by level~-- the resulting tree~$T'$ will be the same.%
   \COMMENT{}

\medbreak

If $\F$ is standard, the $\F$-tangle structure trees we considered in \cref{sec:TSTexistence} were all efficient, because they were thoroughly ordered (\cref{lem:E}). The irreducible structure trees we obtained from them in \cref{thm:irreducible} may no longer be thoroughly ordered: as noted before, this property can get lost in the contraction process. But the efficiency of the original structure trees is maintained:

\begin{lemma}\label{contractionpreservesefficiency}
    If $T'$ is obtained from an efficient tangle structure tree $T$ as in \cref{thm:irreducible}, then $T'$ too is efficient.
\end{lemma}

\begin{proof}
Consider $T^{i+1} = T^i_{w\to v}$ as in \cref{thm:irreducible} for $i=1,\dots,n-1$ in turn. Assuming that $T^i$ is efficient, we have to show that so is~$T^{i+1}$. Write $\beta^{i+1} := \beta^i_{w\to v}$, where $\beta^i$ is the edge labelling of~$T^i$. Given any leaf~$\ell$ of~$T^{i+1}$, we have to show that $\beta^{i+1}_\ell$ is efficient in~$\ucl(\beta^{i+1}_\ell)$. This follows from the fact that $\ell$ is also a leaf of~$T^i$, the assumed efficiency of~$T^i$, and the fact that $\beta^{i+1}_\ell\subseteq \beta^i_\ell$.%
   \COMMENT{}
\end{proof}

\COMMENT{}

\medskip\section{The $\F$-tangle structure theorems}\label{sec:MainThm}

\noindent
   In this section we summarize what we have shown so far, in a few, concise, statements that we think of as the main results of this paper.

The first of these tells us when $\F$-tangle structure trees exist:

\begin{theorem}\label{thm:main}
Let $\vS$ be an ordered separation system. If a set~$\F$ is standard and rich for~$\vS$, then $\vS$ has an efficient and irreducible ordered $\F$-tangle structure tree.
\end{theorem}

\begin{proof}
By \cref{TreeConstruction}, $\vS$~has a thoroughly ordered tangle structure tree $(T,r,\beta)$. By \cref{lem:E} it is efficient. The desired structure tree can be obtained from $(T,r,\beta)$ by \cref{thm:irreducible}; it is efficient by \cref{contractionpreservesefficiency}.
\end{proof}

The tangle structure trees constructed for the proof of \cref{thm:main} have a host of properties designed to help with the tangle analysis of a given separation system. It might have been natural to list these properties as part of \cref{thm:main}. However, since the latter is an existence theorem, this would have allowed us to establish these properties just for the structure tree we constructed in its proof.

We would like to stress the fact that those properties are common to all efficient and irreducible ordered tangle structure trees: these terms were designed to encode precisely those properties. Let us show this next:

\begin{theorem}\label{thm:maindetails}
Let $\vS$ be an ordered separation system. Every efficient ordered%
   \COMMENT{}
   $\F$-tangle structure tree $(T,r,\beta)$ for~$\vS$ has the following properties.
   \begin{enumerate}\itemsep=3pt
   \item For every leaf $\ell$ of~$T$, the set~$\beta_\ell$ is consistent in~$S$. If $\beta_\ell$ has no subset in~$\F$, then $\ucl(\beta_\ell)$ is a tangle of~$S$, and every $\F$-tangle of~$S$ arises in this way.
   \item For every non-leaf node~$v$ of~$T\!$, the set~$\beta_v$ has no subset in~$\F\!$, and $s_v$ is not oriented by~$\ucl(\beta_v)$. Both its orientations $\vs\in\{\svv,\vsv\}$ are minimal in~$\beta_v\cup\{\vs\}$ in the partial ordering of~$\vS$. Its order~$|s|$ is maximal in $\{\,|s_u| : u\le_r v\}$.%
   \COMMENT{}
   \item For every orientation~$\tau$ of~$S$ there is a unique leaf $\ell$ of~$T$ such that $\beta_\ell\subseteq\tau$. If $\tau$ is consistent, then $\ucl(\beta_\ell)\subseteq\tau$ and $\beta_\ell$ is efficient in~$\tau$.
   \item For every $\F$-tangle~$\tau$ of~$S$ there is a unique leaf $\ell$ of~$T$ such that $\beta_\ell\subseteq\tau$. Then $\tau = \ucl(\beta_\ell)$, and $\beta_\ell$ is efficient in the tangle~$\tau$.
   \item For every consistent orientation~$\tau$ of~$S$ that is not a tangle, the set~$\beta_\ell$ in {\rm(iii)} has a subset~$\sigma$ in~$\F$. Every such~$\sigma$ is efficient in~$\tau$.
   \end{enumerate}
\end{theorem}

\begin{proof}
(i) The first two statements are immediate from \cref{def:TST}: tangle structure trees are consistent, and leaves that are not forbidden are tangle leaves. The third assertion is part of (iv) and will be proved there.

(ii) The first assertion is again part of \cref{def:TST}. The second is \cref{cor:B}. The minimality assertion is \cref{lem:notyetnamedlemma}. The maximality assertion holds, because $(T,r,\beta)$ is ordered.

(iii) By \cref{lem:orichoosesleaf}, $T$~has a unique leaf~$\ell$ such that $\beta_\ell\subseteq\tau$. If $\tau$ is consistent, then $\ucl(\beta_\ell)\subseteq \tau$ by \cref{lem:closure}. The claim now follows by \cref{lem:thorough}\,\ref{efficient_and_ordering_2}.

(iv) We have to show that the leaf~$\ell$ from~(iii) satisfies the inclusion $\ucl(\beta_\ell)\subseteq\tau$ provided there with equality. As $\tau$ has no subset in~$\F$, the leaf~$\ell$ is not forbidden. It is therefore a tangle leaf, which means that $\ucl(\beta_\ell)$ is a tangle of~$S$. Since distinct tangles of~$S$ cannot contain one another, this tangle can only be~$\tau$.%
   \COMMENT{}

(v) Consider the leaf~$\ell$ provided for~$\tau$ by~(iii). If $\tau$ is not a tangle, then neither is~$\ucl(\beta_\ell)\subseteq\tau$: if it was, then $\ucl(\beta_\ell)$ and~$\tau$ would both be orientations of~$S$, implying $\ucl(\beta_\ell)=\tau$ with a contradiction.%
   \COMMENT{}
   Hence $\beta_\ell$ has a subset in~$\F$, by~(i). The efficiency claim follows from~(iii), since all subsets of~$\beta_\ell$ are efficient in~$\tau$ if~$\beta_\ell$ is.
\end{proof}

Let us briefly address the question of how efficient our structure trees are in displaying the $\F$-tangles of~$S$, as well as certificates from~$\F$ for consistent orientations of~$S$ that are not $\F$-tangles.

There is a unique way to encode a single tangle~$\tau$ efficiently: we have to know its minimal elements, but only these, since the orientations of all the other separations in~$S$ are then determined by consistency of~$\tau$.%
   \footnote{Apply \cref{lem:closure-consistent} to the set~$\sigma$ of its minimal elements to reobtain the entire tangle as~$\tau=\ucl(\sigma)$.}%
   \COMMENT{}
   The structure tree from \cref{thm:main} does not achieve this for all the tangles of~$S$ simultaneously. Indeed, one only needs two separations in~$S$ to see that no structure tree can in general display only the minimal elements of all its tangles.%
   \footnote{For example, let~$S$ consist of just two nested separations, $r$ and~$s$. Let $\F=\emptyset$, so that all three consistent orientations of~$S$ are $\F$-tangles. The orientations $\vr$ of~$r$ and $\sv$ of~$s$ that point to each other then form a tangle in which they are both minimal, so they must both label an edge of~$T$. Now if we flip either one of them, we obtain a tangle in which the other is still an element, but no longer minimal. If $|r| < |s|$, say, then $\{\vr,\vs\}$ would be such a tangle which our structure tree displays by using both~$\vr$ and~$\vs$ as labels of a path from the root to a tangle leaf, although~$\vs$ alone already determines this tangle. The tangle $\{\rv,\sv\}$ would be displayed by the same tree more efficiently, with $\rv$ labelling an edge from the root to a leaf~$\ell$, with $\beta_\ell = \{\rv\}$ and $\ucl(\beta_\ell) = \{\rv,\sv\}$.}

However, since our structure tree~$T$ in \cref{thm:main} is irreducible,%
   \COMMENT{}
   it does the next best thing. By \cref{thm:maindetails}\,(iv), the minimal elements of every tangle~$\tau$ of~$S$ are displayed as labels along the path~$rT\ell$, where $\ell$ is the tangle leaf with $\tau = \ucl(\beta_\ell)$. Conversely, every label $\vsv$ of an edge of~$T$ is minimal in {\em some\/} tangle of~$S$ (with a tangle leaf~$\ell > v$), or indispensable as an element of the subsets $\sigma\in\F$ of~$\beta_\ell$ that certify that $\ucl(\beta_\ell)$ is not a tangle, for some forbidden leaf~$\ell > v$.

\medbreak

\cref{thm:main} is best possible in another sense too. Its characterizing condition for the existence of tangle structure trees is, essentially, that $\F$ must be rich. This cannot be weakened: if our order function on~$S$ is injective, which in real-world applications is the rule rather than the exception, then by \cref{lemma-twoooo} the existence of any ordered tangle structure tree (efficient and irreducible or not) implies that $\F$ is rich. The richness condition on~$\F$, thus, captures precisely what we need for a tangle structure tree to exist.

If~$\F$ is rich not only for~$\vS$ but also for some (or all)~$\vSk$, \cref{thm:main} applied in~$\vSk$ says that there exist efficient and irreducible ordered tangle structures trees of all these~$\vSk$. But we do not have to obtain these by independent applications of \cref{thm:main}, separately for each~$k$: we can find them all in the original tangle structure tree~$T$ for~$\vS$ provided by \cref{TreeConstruction}, before we applied \cref{thm:irreducible} to `reduce' it. Let us now see how.

Given an ordered separation tree $(T,r,\beta)$ on~$\vS$ and $k\in\R$, let
 $$(T,r,\beta)|_k = (T|_k,r,\beta|_k)$$%
   \COMMENT{}
   be defined as follows. Since $T$ is ordered, its edges~$e$ with $|\beta(e)| < k$ form a subtree of~$T$ rooted at~$r$, which we take as~$T|_k$. As its edge labelling~$\beta|_k$ we take the restriction of~$\beta$ to these edges.

\begin{theorem}\label{thm:Sk}
Let $(T,r,\beta)$ be a thoroughly ordered $\F$-tangle structure tree of~$\vS$. Let $k\in\R$, and assume that $\F$~is rich for~$\vSk$.%
   \COMMENT{}
   Then $(T,r,\beta)|_k$ is a thoroughly ordered $\F$-tangle structure tree~of~$\vSk$.%
   \COMMENT{}
\end{theorem}

\begin{proof}
It is clear from the definition of~$T|_k$ that it inherits from~$T$ the properties of being a thoroughly ordered and consistent separation tree on~$\vSk$.%
   \COMMENT{}
   Moreover, non-leaf nodes~$v$ of~$T|_k$ are non-leaf nodes also of~$T$, so $\beta_v$ has no subset in~$\F$.%
   \COMMENT{}

For a proof that $T|_k$ is a tangle structure tree it remains to show that every leaf~$v$ of~$T|_k$ is either a tangle leaf or a forbidden leaf of~$T|_k$. If it is neither, then by \cref{coreidea}%
   \COMMENT{}
   applied to the closure~$\ucl(\beta_v)_k$ of~$\beta_v$ in~$\vSk$ does not orient all of~$S_k$.

Such a node~$v$ cannot be a leaf of~$T$. Indeed, since $\beta_v$ has no subset in~$\F$, it would be a tangle leaf of~$T$, so $\ucl(\beta_v)$ would be a tangle of~$S$, and $\ucl(\beta_v)\cap\vSk = \ucl(\beta_v)_k$ would be a tangle of~$S_k$ (contrary to our assumption).%
   \COMMENT{}
   So there exists an $s_v\in S$ whose orientations label the edges $e$ and~$e'$ from~$v$ to its children in~$T$.

Since $T$ is thoroughly ordered and there exists a separation in~$S_k$ not oriented by~$\ucl(\beta_v)_k = \ucl(\beta_v)\cap \vSk$, we know that~$|s_v| < k$. But this means that $e$ and~$e'$ are edges of~$T|_k$, contradicting the fact that $v$ is a leaf of~$T|_k$.
\end{proof}

Note that the structure trees~$T|_k$ found inside~$T$ by \cref{thm:Sk} are nested: for $i < j$ we have $T|_i\subseteq T|_j$, indeed $T|_i = (T|_j)|_i$. This hierarchy of structure trees commutes with the following well-known hierarchy of tangles. For all integers $i<j$, every tangle~$\tau_j$ of~$S_j$ induces a tangle~$\tau_i$ of~$S_i$, as $\tau_i = \tau_j\cap\vSi$. A~tangle of any~$S_k$ that is not induced by a tangle of any $S_{k'}\supsetneq S_k$%
   \COMMENT{}
   is {\em maximal\/} in~$\vS$. Every tangle in~$\vS$,%
   \COMMENT{}
   say~$\tau_k$ of order~$k$, now corresponds to a tangle leaf~$\ell(\tau_k)$ of~$T|_k$ in that $\tau_k = \ucl(\beta_{\ell(\tau_k)})_k$. All these~$\ell(\tau_k)$ are nodes of~$T$. 

Conversely, if our order function on~$S$ is injective, then {\em every\/} non-leaf node~$v$ of~$T$ is such a tangle leaf~$\ell(\tau_k)$ of~$T|_k$ for some tangle~$\tau_k$ in~$\vS$, since $\beta_v$ has no subset in~$\F$ by \cref{thm:maindetails}\,(ii): simply choose $k$ just big enough that $\beta(uv) < k$ for $u$ the parent of~$v$.%
   \COMMENT{}
   The positions of these nodes in~$T$ then commute with the relationship of their tangles: if $\tau_j$ induces~$\tau_i$ for $i\le j$, then $\ell(\tau_i)\le_r \ell(\tau_j)$ in the tree-order of~$T$.

Thus, all the maximal tangles in~$\vS$ are displayed by the same structure tree:

\begin{corollary}\label{cor:tanglesinS}
Let $(T,r,\beta)$ be a thoroughly ordered $\F$-tangle structure tree of~$\vS$. Assume that the order function on~$S$ is injective, and that $\F$~is rich for every~$\vSk$.

Let $T'\!$ be obtained from~$T\!$ by deleting any pairs $\ell,\ell'$ of forbidden leaves of~$T\!$ that are children of the same node.%
   \COMMENT{}
   Then $T'\!$~is a thoroughly ordered consistent separation tree on~$\vS$,%
   \COMMENT{}
   which displays precisely the maximal%
   \COMMENT{}
   $\F$-tangles in~$\vS$ as~$\ucl(\beta_\ell)_k$%
   \COMMENT{}
   for leaves~$\ell$ of~$T'\!$ and suitable~$k$.\qed
\end{corollary}

\noindent
   See~\cite[Section~6]{TSTs2} for more details on~$T'$ and applications of \cref{cor:tanglesinS}.

\medbreak

We remark that it is not enough in \cref{thm:Sk} and \cref{cor:tanglesinS} to assume that $\F$ is rich for~$\vS$, rather than for all the relevant~$\vSk$. To see this, recall the example with $\F = \RR$ discussed after \cref{thm:generalduality}. Let $k\in\R$ be large enough that $S=S_k$ for the separation system~$\vS$ discussed there. The unique thoroughly ordered separation tree~$T$ constructed there was not an $\F$-tangle structure tree, because it had leaves that were neither forbidden nor tangle leaves. Correspondingly, $\F$~was not rich for this~$\vS=\vSk$.

Now extend this separation system by adding a new separation~$x$ of order~$k$, making both its orientations incomparable with all the elements of~$\vS$, and adding them to~$\F$ as singleton sets $\{\vx\}$ and~$\{\xv\}$. This makes the extended~$\F'$ rich for the extended~$\vSdash$, because $\{\vx\}$ and~$\{\xv\}$ are efficient and every orientation of~$S'$ contains one of them. The unique $\F'$-tangle structure tree of~$\vSdash$ is obtained from~$T$ by adding two new leaves at every leaf of~$T$, labelling their incident edges with~$\vx$ and~$\xv$. Note that all these leaves are forbidden.

Deleting them, however, as in \cref{cor:tanglesinS}, returns the original~$T$. Contrary to what the corollary claims, $T$~fails to display the maximal $\F'$-tangles in~$\vSdash$ as~$\ucl(\beta_\ell)$ for its leaves~$\ell$, since these are the $\F$-tangles of~$S$.

\medbreak

In the proof of \cref{thm:Sk} we made use of the assumption that $T$ is thoroughly ordered. This holds for the tangle structure trees provided by \cref{TreeConstruction}, and these are already efficient by \cref{lem:E}. But while their efficiency is maintained when we apply \cref{thm:irreducible} to `reduce' them, as in the proof of~\cref{thm:main}, the property of being thoroughly ordered is lost. For the trees~$T|_k$ in \cref{thm:Sk} to be tangle structure trees, however, it is essential that $T$ is thoroughly ordered.%
   \footnote{In the reduction process we contracted or deleted edges whose label was `unnecessary' for $T$ to display the $\F$-tangles of~$S$ and any certificates in~$\F$ of their non-existence. But such edges may have been necessary for $T|_k$ to display the tangles of~$S_k$. As an extreme example, take $\F=\emptyset$ and assume that all the oriented separations that are $\le$-minimal in tangles of~$S$ have order~$k$. While $S$ may have elements of lower order that can be oriented to form lower-order tangles, those tangles will not be visible in~$T$, since $\beta(E(T))\cap\vSk = \emptyset$.}%
   \COMMENT{}

But the tangle structure trees of~$\vSk$ that \cref{thm:Sk} finds in~$T$ can be reduced afterwards: we simply apply \cref{thm:irreducible} to those~$T|_k$. This yields irreducible, efficient,%
   \COMMENT{}
   ordered tangle structure trees for all the~$\vSk$, the same that \cref{thm:main} would give us if we applied it directly in~$\vSk$.

Similarly, we can mimic the reduction process established in \cref{sec:efficient} to contract the separation tree~$T'$ from \cref{cor:tanglesinS}, which displays all the maximal tangles in~$\vS$ as~$\ucl(\beta_\ell)$ for suitable leaves~$\ell$, to an efficient and `irreducible' such tree. All it needs to facilitate this is that we adapt the definition of `tangle leaves', and those of ensuing terms such as `necessary' and `irreducible', to refer to the maximal tangles in~$\vS$, rather than just the tangles of all of~$S$.%
   \footnote{\lineskiplimit -3pt The only (small) addition in substance%
   \COMMENT{}
   needed here occurs early in the proof of \cref{lem:necessity}, where we now have to show also that $\ucl(\beta_\ell\setminus\{\vs\})$ cannot be a maximal tangle in~$\vS$ other than~$\tau$. But this is clear, since $\ucl(\beta_\ell\setminus\{\vs\})\subseteq\tau$.}

\medbreak

Finally, let us apply \cref{thm:main} to the special case that $S$ has no $\F$-tangles. In this case it is desirable to be able to certify this efficiently. 

In the special case that $\F$ consists of stars (of oriented separations; see \cref{sec:defs}), such certificates are known: they are the {\em $S$-trees over~$\F$} from~\cite{TangleTreeAbstract}, which generalize the branch-decompositions introduced by Robertson and Seymour~\cite{GMX} for graph tangles to $\F$-tangles in abstract separation systems.

Our structure trees from \cref{thm:main} provide maximally efficient%
   \COMMENT{}
   certificates for the non-existence of $\F$-tangles for general~$\F$, not necessarily consisting of stars, when all their leaves are forbidden. Let us call such $\F$-tangle structure trees {\em $\F$-trees\/}.\looseness=-1

\begin{theorem}\label{thm:generalduality2}
Let $\vS$ be an ordered%
   \COMMENT{}
   separation system, and let~$\F$ be any set that is standard and rich for~$\vS$. Then exactly one of the following assertions holds:
\begin{enumerate}\itemsep2pt\vskip2pt
  \item\label{item:gd1} there exists an $\F$-tangle of~$S$;
  \item\label{item:gd2} there exists an $\F$-tree of~$\vS$.
\end{enumerate}
In the case of {\rm\labelcref{item:gd2}}, the $\F$-tree can be chosen to be irreducible, efficient and ordered.
\end{theorem}

\begin{proof}
By \cref{thm:C} \labelcref{item:C2}, any $\F$-tree of~$S$ precludes the existence of an $\F$-tangle of~$\vS$, so we cannot have both \labelcref{item:gd1} and \labelcref{item:gd2}.

Let us prove that \labelcref{item:gd1} or \labelcref{item:gd2} holds. By \cref{thm:main}, $\vS$~has an efficient and irreducible ordered $\F$-tangle structure tree~$T$. If \labelcref{item:gd1} fails, then $T$ has no tangle leaf. Then all its leaves are forbidden, so $T$ is an $\F$-tree.
\end{proof}

Let us close this section with an interesting fact exhibited by $\F$-trees. If $S$ has no tangle, then orientations of~$S$ cannot fail to be a tangle just because they are inconsistent: even inconsistent orientations of~$S$ must have a subset in~$\F$. This is clearly not the case when $S$ does have tangles, e.g.\ if $\F=\emptyset$, as soon it contains two nested separations (which we can orient away from each other).

\begin{corollary}
Let $\vS$ be an ordered%
   \COMMENT{}
   separation system, and let~$\F$ be any set that is standard and rich for~$\vS$. If $S$ has no $\F$-tangle, then every orientation of~$S$ has a subset in~$\F$.
\end{corollary}

\begin{proof}
Let $\tau$ be any orientation of~$S$. By \cref{thm:generalduality2}, $S$~has an $\F$-tree $(T,r,\beta)$. Let $rT\ell$ be the maximal path in~$T$ from the root such that $\beta(uv)\in\tau$ for its adjacent vertices ${u<v}$. Since $\ell$ is a forbidden leaf, $\F$~has an element $\sigma\subseteq\beta_\ell\subseteq\tau$.
\end{proof}

\medskip\section{Applications: blocks, profiles, and cluster tangles}\label{sec:apps}

\noindent
   Graph tangles can be expressed as $\F$-tangles in the universe of all separations of a graph~~$G$: take as $\F$ the set~$\T$ of triples $\{\,(A_i,B_i)\mid i=1,2,3\,\}$ of oriented separations such that $G[A_1]\cup G[A_2]\cup G[A_3] = G$. These are not stars. But it is not hard to show~\cite{AbstractTangles} that the $\T$-tangles of a graph are precisely its $\T^*$-tangles, where $\T^*$ consists of the sets in~$\T$ that are stars.%
   \COMMENT{}

Due to this history, $\F$-tangles have so far been studied mostly when $\F$ consisted of stars of separations. Such $\F$-tangles already allow for vast generalizations of the original tangles of graphs. However, there are contexts in which $\F$-tangles occur naturally and $\F$ does not consist of stars. Three of these are:

\smallskip
\begin{itemize}[itemsep=3pt]
\item blocks in graphs, more precisely: $k$-blocks for any $k\in\N$;
\item profiles: the most comprehensive generalization of graph tangles that still admits a tree-of-tangles structure theorem;
\item cluster tangles in large datasets.
\end{itemize}
\smallskip

\noindent
   We shall apply our results to all these in this section, giving a minimum of background for context in each case.

Common to all tangle contexts are  two fundamental theorems about tangles in graphs, which both go back to the original paper of Robertson and Seymour~\cite{GMX} in which tangles were first introduced. These have become known~\cite{DiestelBook25}~as

\smallskip
\begin{itemize}[itemsep=3pt]
\item the tree-of-tangles theorem; and
\item the tangle-tree duality theorem.
\end{itemize}
\smallskip

The tree-of-tangles theorem, in its simplest form, says that all the maximal tangles in a graph can be distinguished%
   \footnote{Recall that a separation {\em distinguishes\/} two tangles if they orient it differently.}
   by a nested set~$T$ of separations, one that can be represented as the set of separations associated with a \td\ of the graph. The separations in~$T$ can be chosen so that they have minimum order for any separation that distinguishes the pair of tangles they distinguish,%
   \COMMENT{}
   and {\em canonically\/}, which means that $T$ is invariant under all the automorphisms of the graph~\cite{DiestelBook25}. These enhanced tree-of-tangle theorems have been further generalized to abstract separation systems~\cite{ProfilesNew,AbstractTangles,carmesin2022entanglements}.

The tangle-tree duality theorem, in its simplest form, says that if a graph has no tangle of some given order~$k$ then it has a nested set~$T$ of separations of order~$<k$ that witnesses the non-existence of those tangles: any orientation~$\tau$ of the graph's separations of order~$<k$ will in particular orient~$T$, and thereby the edges of the decomposition tree associated with~$T$, and following those oriented edges will take us to a node whose incident edges correspond to a set in~$\T^*$; thus, $\tau$~is not a tangle. The tangle-tree duality theorem, too, has been generalized to abstract separation systems~\cite{TangleTreeAbstract,AbstractTangles}. There are also unifications of both theorems into one~\cite{JoshRefining,SARefiningInessParts,SARefiningEssParts}.

The certificates for the non-existence of an $\F$-tangle offered by the tangle-tree duality theorems require that $\F$ consist of stars. Indeed, since $T$ is nested, the edges at any node of the decomposition tree that represents~$T$ will map to a nested set of separations too, and orienting them towards that node will yield a star of separations. Hence if we look for certificates for the non-existence of $\F$-tangles when $\F$ does not consist of stars, we need a fresh start. Our \cref{thm:generalduality2} offers such an alternative.

\medbreak

For our first application, let us now look at blocks. For any $k\in\N$, a {\em $k$-block\/} in a graph is any maximal set of at least~$k$ vertices no two of which can be separated by $<k$ vertices. Every $k$-block~$X$ in a graph~$G$ defines an $\FBk$-tangle of~$S_k$ for
 $$\FBk := \left\{\,\sigma\subseteq\vS : \Big|\bigcap\{\,B\mid (A,B)\in\sigma\}\Big| < k\,\right\},$$
 where $S$ is the set of all separations of~$G$: just orient every separation in~$S_k$ towards~$X$. Conversely, given any $\FBk$-tangle~$\tau$ of~$S_k$ in~$G$, the set $\bigcap\{\,B\mid (A,B)\in\tau\}$ is a $k$-block in~$G$, which in turn defines~$\tau$ as indicated.

One of the earliest applications of what later became $\F$-tangles was that finite graphs have canonical \td s that distinguish all their $k$-blocks~\cite{confing}.%
   \COMMENT{}
   This extended a classical  decomposition theorem of Tutte for 2-blocks. But there was no tangle-tree duality theorem for blocks, since $\FBk$ does not consist of stars.

It was shown much later that blocks can in fact be captured by $\F$-tangles with $\F$ consisting of stars~\cite{EberenzMaster,ProfileDuality,AbstractTangles}. But these~$\F$ have to be constructed recursively for each given graph, and are therefore highly artificial: tailor-made for the purpose, but not defining $\F$-tangles that can be understood from their definition. As a consequence, the certificates they yield in theory for the non-existence of blocks are of little use in practice.
Our \cref{thm:generalduality2}, however, offers very tangible certificates:

\begin{theorem}\label{thm:blocks}
Let $S$ be the set of separations of a graph~$G$, with the usual order function~\cite{DiestelBook25}, and let $k\in\N$. Then exactly one of the following assertions holds:
\begin{enumerate}\itemsep2pt\vskip2pt
  \item\label{item:gd1} $G$ has a $k$-block;
  \item\label{item:gd2} there exists an $\FBk$-tree of~$\vSk$.
\end{enumerate}
In the case of {\rm\labelcref{item:gd2}}, the $\FBk$-tree can be chosen to be irreducible, efficient, and ordered.
\end{theorem}

\begin{proof}
We have already seen that $k$-blocks induce $\FBk$-tangles of~$S_k$ and vice versa. We can thus replace~\labelcref{item:gd1} with the assertion that $S_k$ has a $\FBk$-tangle, and proceed to apply \cref{thm:generalduality2}.

Let us check the theorem's premise. The set $\FBk$ is standard for~$\vSk$, since the all co-trivial oriented separations in~$\vSk$ are of the form~$(V,A)$, where $V=V(G)$ and $|A|<k$.%
   \footnote{In addition, the set~$A$ has to lie in the separator of another separation of order~$<k$.}
   Clearly, then, $\{\{(V,A)\}\}\in\FBk$.

For a proof that $\FBk$ is rich for~$\vSk$, it suffices by \cref{minimization} to show that $\FBk$ is closed under minimization in~$\vS$. Given any $\sigma\in\F$, let $\sigma'$ be obtained from~$\sigma$ by replacing every $(A,B)\in\sigma$ by some $(A',B')\le (A,B)$ from~$\vS$.%
   \footnote{Recall that, for graph separations, $(A',B')\le (A,B)$ if and only if $A'\supseteq A$ and $B'\subseteq B$.}%
   \COMMENT{}
   As 
 $$\bigcap\,\big\{\,B' : (A',B')\in\sigma'\big\}\subseteq \bigcap\,\big\{\,B : (A,B)\in\sigma\,\big\},$$
  our assumption of $\sigma\in\FBk$ clearly implies that $\sigma'\in\FBk$, as required.
\end{proof}

Next, profiles. Profiles of graphs were introduced by Hundertmark~\cite{profiles}, who noticed that the above-mentioned `tree-of-tangles' theorem for blocks~\cite{confing} required in its proof much less information about blocks than is provided in their definition. All that is needed is that the orientation of~$S_k$ which a $k$-block induces%
   \footnote{orient every separation towards the side that contains the given $k$-block}
   avoids the set~$\F$ of triples of the form $\{\vr,\vs,\rv\lor\sv\}$. In the context of graphs, this supremum is taken in the universe of all the separations of the given graph. 

More generally, given any separation system~$\vS$ in some universe~$\vU$ of separations, we let
 $$\P = \big\{\,\{\vr,\vs,\rv\lor\sv\} : \vr,\vs\in\vU\,\big\},$$%
   \COMMENT{}
 and refer to the $\P$-tangles of~$S$ as its {\em profiles\/}. A~profile is {\em regular\/} if none of its elements is small.

Hundertmark's discovery was seminal: since graph tangles are profiles too, the proof of the canonical tree-of-tangles theorem for blocks~\cite{confing}, once rewritten for profiles~\cite{ProfilesNew}, yielded {\em canonical\/} trees of graph tangles~-- something that Robertson and Seymour~\cite{GMX} had been unable to achieve directly by their tangle methods.

So what about tangle-tree duality for profiles? It was the hope to find such a dichotomy theorem for profiles that motivated the original generalization of graph tangles to $\F$-tangles in abstract separation systems~\cite{TangleTreeAbstract}, since profiles had this form by definition. But the hope proved elusive. We did indeed find a duality theorem for $\F$-tangles in~\cite{TangleTreeAbstract}, but this once more required that $\F$ consisted of stars. This was good enough for tangles, since $\T$-tangles coincide with $\T^*$-tangles, as well as for some further duality theorems in graphs and matroids~\cite{TangleTreeGraphsMatroids}, but there was no tangle-tree dichotomy theorem for profiles among these.

As in the case of blocks it is possible to construct, recursively, sets~$\F$ of stars in submodular separation systems embedded in a%
   \COMMENT{}
   universe so that its profiles are precisely its~$\F$-tangles. One can then apply~\cite{TangleTreeAbstract} to obtain a tangle-tree duality theorem for profiles, encoded as such $\F$-tangles~\cite{ProfileDuality}.

Alternatively, we can now obtain certificates for the non-existence of profiles directly from \cref{thm:generalduality2}.

As preparation, we need a lemma. Given a separation system~$\vS$ in some universe~$\vU$, we call a profile of~$S$ {\em strong\/} if it has no subset in
  $$\FP = \big\{\,\{\vr,\vs,\tv\} : \vr,\vs,\vt\in\vU\hbox{ and }\tv\le\rv\lor\sv\,\big\}.$$

Note that strong profiles are regular: if $\vr\le \rv$ then no strong profile contains~$\vr$, since $\vr\le\rv\lor\rv$ puts $\{\vr\}$ in~$\FP$.
If $\vS$ is submodular and $\vU$ is distributive,%
   \COMMENT{}
   the converse holds too:

\begin{lemma}\label{strongprofiles}
Every regular profile of a submodular separation system in a distributive universe of separations is strong.%
   \COMMENT{}
\end{lemma}

\begin{proof}
Let $(\vS,\le)$ be any submodular separation system in a distributive lattice~$\vU$. Suppose $S$ has a regular profile~$\tau$ that is not strong. Then there exist $\vr, \vs, \tv$ in~$\tau$ such that $\tv \leq \rv \lor \sv$. Choose this triple with $\tv$ as large as possible under~$\le$.

Our aim will be to show that both $\vt\lor\vr$ and $\vt\lor\vs$ lie in~$\vS$.%
   \COMMENT{}
   For if they do, they must lie in~$\tau$ by the consistency of~$\tau$, since $\vt\lor\vr\ge\vr\in\tau$ and $\vt\lor\vs\ge\vs\in\tau$ and $\tau$ is regular.%
   \footnote{Compare \cref{lem:closure} and the remark preceding it.}
   Then, as
 $$(\vt \lor \vs)^* \lor (\vt \lor \vr)^* = (\tv \land \sv) \lor (\tv \land \rv) = \tv \land (\sv \lor \rv) = \tv$$%
   \COMMENT{}
by distributivity, we shall have $\tau\supseteq \{\vt \lor \vr, \vt \lor \vs, \tv\} \in \P$, which contradicts our assumption that $\tau$ is a profile. 

So let us show that $\vt\lor\vr\in\vS$; the proof of $\vt\lor\vs\in\vS$ is analogous. By submodularity we have $\vt\lor\vr\in\vS$ as desired unless $\vt\land\vr\in\vS$, so let us assume this. Then also $\vS\owns (\vt\land\vr)^* = \tv\lor\rv \ge \tv\in\tau$, so $\tv\lor\rv\in\tau$ by the consistency and regularity$^{\thefootnote}$ of~$\tau$.

If $\tv\lor\rv = \tv$,%
   \COMMENT{}
  then $\tv\ge\rv$%
  \COMMENT{}
  and thus $\vt\le\vr$, giving $\vt\lor\vr = \vr\in\vS$ as desired. But otherwise $\tv\lor\rv > \tv\in\tau$, which contradicts the maximality of~$\tv$ in our choice of $\vr,\vs,\tv\in\tau$ with $\tv\le\rv\lor\sv$, since also $\tv\lor\rv\le\rv\lor\sv$ by definition of the supremum~$\tv\lor\rv$.%
   \COMMENT{}
   \end{proof}

\cref{strongprofiles} implies that the regular profiles of a submodular separation system in a distibutive universe are precisely its~$\FP$-tangles. Indeed, being profiles, they are~$\P$-tangles and therefore consistent, and they avoid~$\FP$ since, by the lemma, they are strong, so they are $\FP$-tangles. Conversely, every $\FP$-tangle avoids~$\P\subseteq\FP$ and is therefore a profile, and it is strong since $\FP$-tangles cannot contain small separations, as noted earlier.

\begin{theorem}\label{thm:profiles}
Let $\vS$ be any submodular separation system in a distibutive universe of separations. Then exactly one of the following assertions holds:
\begin{enumerate}\itemsep2pt\vskip2pt
  \item there exists a regular profile of~$S$;
  \item there exists an $\FP$-tree of~$\vS$.
\end{enumerate}
In the case of {\rm\labelcref{item:gd2}}, the $\FP$-tree can be chosen to be irreducible, efficient, and ordered with respect to any order function on~$S$.
\end{theorem}

\begin{proof}
Our aim is us apply \cref{thm:generalduality2} to~$\vS$ with $\F=\FP$, and richness referring to any order function on~$S$ we may choose.%
   \COMMENT{}
   As remarked after \cref{strongprofiles}, the regular profiles of~$S$ are precisely its $\FP$-tangles. It remains to check that $\FP$ is standard and rich for~$\vS$.

For a proof that $\FP$ is standard, we have to show that $\{\sv\}\in\FP$ whenever $\vs$ is trivial in~$\vS$. Let $r\in S$ witness that $\vs$ is trivial. Then $\vr,\rv < \vs$, so $\sv < \vr < \vs$. Hence $\sv\le \vs\lor\vs$, which puts $\{\sv\}$ in~$\FP$.

For a proof that $\FP$ is rich for~$\vS$, it suffices by \cref{minimization} to show that $\FP$ is closed under minimization in~$\vS$. Given any set $\sigma = \{\vr,\vs,\tv\}\in\FP$, consider any set $\sigma' = \{\vrdash,\vsdash,\tvdash\}$ with $\vrdash\le\vr$ and $\vsdash\le\vs$ and $\tvdash\le\tv$. Then $\sigma'\in\FP$ as desired, since $\tvdash\le\tv\le\rv\lor\sv\le\rvdash\lor\svdash$; the middle inequality comes from $\sigma\in\FP$.%
   \COMMENT{}
  \end{proof}

Finally, let us apply our results to tangles defined by clusters in large datasets. The basic setup is that $\vS$ is the set~$2^V\!$ of all subsets of our data set~$V\!$, the involution~* given by complementation in~$V\!$. Then $S$~is the set of bipartitions $s = \{A,B\}$ of~$V\!$, whose orientations are $\vs=A$ and~$\sv=B$ (say). We then choose an order function on~$S$ that assigns high values to separations that divide many pairs of elements of~$V\!$ which we consider as `close',%
   \footnote{Choosing this order function to fit the application context allows for considerable variety in using tangles for clustering~\cite[Chapter~9]{TangleBook}. See~\cite{Big5clustering} for a concrete example.}
   the idea being that separations of low order cannot divide many close pairs and therefore cannot cut through a cluster, only chip off a few points. This, then, implies that large clusters induce $\F$-tangles for~$\F$ defined as\looseness=-1
 $$\FCn := \big\{\,\{\vr,\vs,\vt\}: |\vr\cap\vs\cap\vt| < n\,\big\},$$
 where the {\em agreement\/} value~$n$ is chosen to fit the context.

Thus, $\FCn$-tangles in datasets are a bit like $\FBk$-tangles in graphs, except that $n$ is now a value separate from the order~$k$ of the tangles considered. An important difference is that, unlike $k$-blocks in graphs, the clusters captured by those tangles are not normally subsets~$X$ of~$V\!$ that we can, or even wish to, specify explicitly: they can be `fuzzy', and do not have to lie entirely on the side of the separations chosen by the tangle they induce, only mostly. Quantitatively, while the intersection of any three elements of an $\FCn$-tangle contains at least~$n$ datapoints, the intersection of all of them will likely be empty~-- unlike in the case of $k$-blocks in graphs, which are equal to the intersection of {\em all\/} the sides of the graph's separations chosen by the tangle.\looseness=-1

In other contexts than generic clustering one sometimes looks for $\FCn$-tangles not of the set~$S$ of all the bipartitions of~$V\!$, but of some specific, hand-designed set~$S$ of bipartitions of~$V\!$. For example, $S$~might be a questionnaire whose elements partition a population~$V\!$ of individuals that have answered it into those that answered yes and those that answered no. In such a context, $\FCn$-tangles can be interpreted as `typical' ways to answer the questions in~$S$, as {\em mindsets\/} found in the population~$V\!$ regarding the questions in~$S$. See~\cite{TangleBook} for more.

$\!$In all these contexts we can now use \cref{thm:main} to efficiently display the clusters or mindsets in~$V\!$ as determined by~$S$, as long as $\FCn$ is standard and rich for~$\vS$. If we are interested specifically in certifying that our dataset has no clusters of some desired density%
   \COMMENT{}
   at all, we can similarly apply \cref{thm:generalduality2}.

The sets~$\FCn$ are standard for all~$\vS\ne\emptyset$ and $n>0$: then $\{\emptyset\}\in\FCn$, and only~$\emptyset$ can be%
   \COMMENT{}
   co-trivial in~$\vS$. Also, $\FCn$~is clearly closed under minimization in~$\vS$: replacing the elements~$\vs$ of a triple~$\sigma$ with subsets~$\vsdash\le\vs$ leads to $\bigcap\sigma'\subseteq\bigcap\sigma$, so $\sigma\in\FCn$ implies $\sigma'\in\FCn$. By \cref{minimization}, therefore, $\FCn$~is rich for~$\vS$.

\cref{thm:main,thm:generalduality2} thus have the following instances in our context:

\begin{corollary}\label{cor:main}
Let $\vS$ be any non-empty set of subsets of a set~$V\!$ that is closed under taking complements in~$V\!$. Let any order function on~$S$ be given, and let $n>0$. Then $\vS$~has an efficient and irreducible ordered $\FCn$-tangle structure tree.\qed
\end{corollary}

\begin{corollary}\label{cor:generalduality2}
Let $\vS$ be any non-empty set of subsets of a set~$V\!$ that is closed under taking complements in~$V\!$. Let any order function on~$S$ be given, and let $n>0$.\penalty-200\ Then exactly one of the following assertions holds:
\begin{enumerate}\itemsep2pt\vskip2pt
  \item\label{item:gd1} there exists a $\FCn$-tangle of~$S$;
  \item\label{item:gd2} there exists a $\FCn$-tree of~$\vS$.
\end{enumerate}
In case~{\rm\labelcref{item:gd2}}, the $\FCn$-tree can be chosen to be irreducible, efficient, and ordered.\qed
\end{corollary}

\bibliographystyle{abbrv}
\bibliography{collective}

\end{document}

%% file: letterswitharrowsRhd.tex
\def\lowfwd #1#2#3{{\mathop{\kern0pt #1}\limits^{\kern#2pt\raise.#3ex
\vbox to 0pt{\hbox{$\scriptscriptstyle\rightarrow$}\vss}}}}
\def\lowbkwd #1#2#3{{\mathop{\kern0pt #1}\limits^{\kern#2pt\raise.#3ex
\vbox to 0pt{\hbox{$\scriptscriptstyle\leftarrow$}\vss}}}}
\def\lowfbwd #1#2#3{{\mathop{\kern0pt #1}\limits^{\kern#2pt\raise.#3ex
\vbox to 0pt{\hbox{$\scriptscriptstyle\leftrightarrow$}\vss}}}}
\def\fwd #1#2{{\lowfwd{#1}{#2}{15}}}
\def\ve{\kern-1.5pt\lowfwd e{1.5}2\kern-1pt}
\def\vedash{{\mathop{\kern0pt e\lower.5pt\hbox{${}
     \scriptstyle'$}}\limits^{\kern0pt\raise.02ex
     \vbox to 0pt{\hbox{$\scriptscriptstyle\rightarrow$}\vss}}}}
\def\ev{\kern-1pt\lowbkwd e{0.5}2\kern-1pt}
\def\vf{\kern-2pt\lowfwd f{2.5}2\kern-1pt}
\def\vfdash{{\mathop{\kern0pt f\raise 1pt\hbox{${}
     \scriptstyle'$}}\limits^{\kern2pt\raise.02ex
     \vbox to 0pt{\hbox{$\scriptscriptstyle\rightarrow$}\vss}}}}

\def\vr{\lowfwd r{1.5}2}
\def\rv{\lowbkwd r02}

\def\vrdash{{\mathop{\kern0pt r\lower.5pt\hbox{${}
     \scriptstyle'$}}\limits^{\kern0pt\raise.02ex
     \vbox to 0pt{\hbox{$\scriptscriptstyle\rightarrow$}\vss}}}}
\def\rvdash{{\mathop{\kern0pt r\lower.5pt\hbox{${}
     \scriptstyle'$}}\limits^{\kern0pt\raise.02ex
     \vbox to 0pt{\hbox{$\scriptscriptstyle\leftarrow$}\vss}}}}
\def\vrdashp{{\mathop{\kern0pt r_p\kern-4pt\lower.5pt\hbox{${}
     \scriptstyle'$}}\limits^{\kern0pt\raise.02ex
     \vbox to 0pt{\hbox{$\scriptscriptstyle\rightarrow$}\vss}}}\,}
\def\rvdashp{{\mathop{\kern0pt r_p\kern-4pt\lower.5pt\hbox{${}
     \scriptstyle'$}}\limits^{\kern0pt\raise.02ex
     \vbox to 0pt{\hbox{$\scriptscriptstyle\leftarrow$}\vss}}}\,}
\def\vrddash{{\mathop{\kern0pt r\lower.5pt\hbox{${}
     \scriptstyle''$}}\limits^{\kern0pt\raise.02ex
     \vbox to 0pt{\hbox{$\scriptscriptstyle\rightarrow$}\vss}}}}

\def\vs{\hskip-1pt\lowfwd s{1.5}1}
\def\sv{{{\hskip-1pt\lowbkwd s{1}1}\hskip-1pt}}
\def\vsdash{{\mathop{\kern0pt s\lower.5pt\hbox{${}
     \scriptstyle'$}}\limits^{\kern0pt\raise.02ex
     \vbox to 0pt{\hbox{$\scriptscriptstyle\rightarrow$}\vss}}}}
\def\svdash{{\mathop{\kern0pt s\lower.5pt\hbox{${}
     \scriptstyle'$}}\limits^{\kern0pt\raise.02ex
     \vbox to 0pt{\hbox{$\scriptscriptstyle\leftarrow$}\vss}}}}
\def\vsddash{{\mathop{\kern0pt s\lower.5pt\hbox{${}
     \scriptstyle''$}}\limits^{\kern0pt\raise.02ex
     \vbox to 0pt{\hbox{$\scriptscriptstyle\rightarrow$}\vss}}}}
\def\svddash{{\mathop{\kern0pt s\lower.5pt\hbox{${}
     \scriptstyle''$}}\limits^{\kern0pt\raise.02ex
     \vbox to 0pt{\hbox{$\scriptscriptstyle\leftarrow$}\vss}}}}
\def\vsdashp{{\mathop{\kern0pt s_p\kern-4pt\lower.5pt\hbox{${}
     \scriptstyle'$}}\limits^{\kern0pt\raise.02ex
     \vbox to 0pt{\hbox{$\scriptscriptstyle\rightarrow$}\vss}}}\,}
\def\svdashp{{\mathop{\kern0pt s_p\kern-4pt\lower.5pt\hbox{${}
     \scriptstyle'$}}\limits^{\kern0pt\raise.02ex
     \vbox to 0pt{\hbox{$\scriptscriptstyle\leftarrow$}\vss}}}\,}

\def\vsidash{{\mathop{\kern0pt s_i\kern-3.5pt\lower.3pt\hbox{${}
     \scriptstyle'$}}\limits^{\kern0pt\raise.02ex
     \vbox to 0pt{\hbox{$\scriptscriptstyle\rightarrow$}\vss}}}}

\def\vsqdash{{\mathop{\kern0pt s_q\kern-3.5pt\lower.3pt\hbox{${}
     \scriptstyle'$}}\limits^{\kern0pt\raise.02ex
     \vbox to 0pt{\hbox{$\scriptscriptstyle\rightarrow$}\vss}}}}

\def\vtdash{{\mathop{\kern0pt t\lower-.5pt\hbox{${}
     \scriptstyle'$}}\limits^{\kern0pt\raise.1ex
     \vbox to 0pt{\hbox{$\scriptscriptstyle\rightarrow$}\vss}}}}
\def\tvdash{{\mathop{\kern0pt t\lower-.5pt\hbox{${}
     \scriptstyle'$}}\limits^{\kern0pt\raise.1ex
     \vbox to 0pt{\hbox{$\scriptscriptstyle\leftarrow$}\vss}}}}

\def\vsu{\lowfwd {s_u}11}
\def\svu{\lowbkwd {s_u}{-1}2}
\def\vsv{\lowfwd {s_v}11}
\def\svv{\lowbkwd {s_v}{-1}2}
\def\vsw{\lowfwd {s_w}11}
\def\svw{\lowbkwd {s_w}{-1}2}
\def\vS{{\hskip-1pt{\fwd S3}\hskip-1pt}} 
\def\vSi{\lowfwd {S_i}11}
\def\vSk{\lowfwd {S_k}11}
\def\vSk{\lowfwd {S_k}11}

\def\vSstar{{\mathop{\kern0pt S\lower-1pt\hbox{$^*$}}\limits^{\kern2pt
     \vbox to 0pt{\hbox{$\scriptscriptstyle\rightarrow$}\vss}}}}
\def\vSdash{{\mathop{\kern0pt S\lower-1pt\hbox{${}
     \scriptstyle'$}}\limits^{\kern2pt\raise.1ex
     \vbox to 0pt{\hbox{$\scriptscriptstyle\rightarrow$}\vss}}}}

\def\vt{\lowfwd t{1.5}1}
\def\tv{\lowbkwd t{1.5}1}

\def\vU{{\vec U}} 

\def\vx{\hskip-1pt\lowfwd x{1.5}1}
\def\xv{{{\hskip-1pt\lowbkwd x{1}1}\hskip-1pt}}

%% file: TSTsArXiv.bbl
\begin{bibdiv}
\begin{biblist}

\bib{SARefiningEssParts}{misc}{
      author={Albrechtsen, Sandra},
       title={Optimal trees of tangles: refining the essential parts},
         how={ar{X}iv:2304.12078},
        date={2023},
}

\bib{SARefiningInessParts}{article}{
      author={Albrechtsen, Sandra},
       title={Refining trees of tangles in abstract separation systems: {I}nessential parts},
        date={2024; ar{X}iv:2302.01808},
     journal={Combinatorial Theory},
      volume={4(1)},
}

\bib{Big5clustering}{article}{
      author={Bergen, {H.v.}},
      author={Diestel, R.},
       title={Traits and tangles: an analysis of the {Big Five} paradigm by tangle-based clustering},
        date={2025; ar{X}iv:2411.18670},
     journal={J.~Math.\ Psychology},
      volume={125 (102920)},
}

\bib{TSTs2}{misc}{
      author={Bergen, {H.v.}},
      author={Diestel, R.},
       title={Tangle structure trees~{II}: trees of tangles and tangle-tree duality},
        date={arXiv:2603.17756},
}

\bib{CDHH13CanonicalAlg}{article}{
      author={Carmesin, J.},
      author={Diestel, R.},
      author={Hamann, M.},
      author={Hundertmark, F.},
       title={Canonical tree-\penalty-200 decompositions of finite graphs {I.~Existence} and algorithms},
        date={2016},
     journal={J. Comb. Theory Ser. B},
      volume={116},
       pages={1\ndash 24; arXiv:1305.4668},
}

\bib{CDHH13CanonicalParts}{article}{
      author={Carmesin, J.},
      author={Diestel, R.},
      author={Hamann, M.},
      author={Hundertmark, F.},
       title={Canonical tree-\penalty-200 decompositions of finite graphs {II.~Essential} parts},
        date={2016},
     journal={J. Comb. Theory Ser. B},
      volume={118},
       pages={268\ndash 283; arXiv:1305.4909},
}

\bib{confing}{article}{
      author={Carmesin, J.},
      author={Diestel, R.},
      author={Hundertmark, F.},
      author={Stein, M.},
       title={Connectivity and tree structure in finite graphs},
        date={2014},
     journal={Combinatorica},
      volume={34},
      number={1},
       pages={1\ndash 35; arXiv:1105.1611},
}

\bib{carmesin2022entanglements}{article}{
      author={Carmesin, J.},
      author={Kurkofka, J.},
       title={Entanglements},
        date={2024},
     journal={J.~Comb.\ Theory, Ser.~B},
      volume={164},
       pages={17\ndash 28; arXiv:2205.11488},
}

\bib{ASS}{article}{
      author={Diestel, R.},
       title={Abstract separation systems},
         how={arXiv:1406.3797},
        date={2018},
     journal={Order},
      volume={35},
       pages={157\ndash 170; arXiv:1406.3797},
}

\bib{TangleBook}{book}{
      author={Diestel, R.},
       title={Tangles: A~structural approach to artificial intelligence in the empirical sciences},
   publisher={Cambridge University Press},
        date={2024. {e-book: {\tt tangles-book.com/book/}}},
}

\bib{DiestelBook25}{book}{
      author={Diestel, R.},
       title={Graph theory \emph{(6th edition)}},
   publisher={Springer-Verlag},
        date={2025},
        note={\\ Electronic edition available at {\small\tt http://diestel-graph-theory.com/}},
}

\bib{ProfileDuality}{article}{
      author={Diestel, R.},
      author={Eberenz, Ph.},
      author={Erde, J.},
       title={Duality theorems for blocks and tangles in graphs},
        date={2017},
     journal={SIAM J.\ Discrete Math.},
      volume={31},
      number={3},
       pages={1514\ndash 1528; ar{X}iv:1605.09139},
}

\bib{AbstractTangles}{article}{
      author={Diestel, R.},
      author={Erde, J.},
      author={Wei{\ss}auer, D.},
       title={Structural submodularity and tangles in abstract separation systems},
        date={2019},
     journal={J. Comb. Theory Ser. A},
      volume={167C},
       pages={155\ndash 180; ar{X}iv:1805.01439},
}

\bib{ProfilesNew}{article}{
      author={Diestel, R.},
      author={Hundertmark, F.},
      author={Lemanczyk, S.},
       title={Profiles of separations: in graphs, matroids, and beyond},
        date={2019},
     journal={Combinatorica},
      volume={39},
      number={1},
       pages={37\ndash 75; arXiv:1110.6207},
}

\bib{TangleTreeGraphsMatroids}{article}{
      author={Diestel, R.},
      author={Oum, S.},
       title={Tangle-tree duality in graphs, matroids and beyond},
        date={2019},
     journal={Combinatorica},
      volume={39},
       pages={879\ndash 910},
        note={arXiv:1701.02651},
}

\bib{TangleTreeAbstract}{article}{
      author={Diestel, R.},
      author={Oum, S.},
       title={Tangle-tree duality in abstract separation systems},
        date={2021},
     journal={Advances in Mathematics},
      volume={377},
       pages={107470; arXiv:1701.02509},
}

\bib{EberenzMaster}{thesis}{
      author={Eberenz, Ph.},
       title={Characteristics of profiles},
        type={Master's Thesis},
         how={MSc dissertation, Hamburg 2015},
        date={Hamburg University, 2015},
}

\bib{ToTfromTTD}{article}{
      author={Elbracht, C.},
      author={Kneip, J.},
      author={Teegen, M.},
       title={Obtaining trees of tangles from tangle-tree duality},
        date={2022},
     journal={J. Combinatorics},
      volume={13},
       pages={251\ndash 287; arXiv:2011.09758},
}

\bib{JoshRefining}{article}{
      author={Erde, Joshua},
       title={Refining a tree-decomposition which distinguishes tangles},
        date={2017},
     journal={SIAM J.\ Discrete Math.},
      volume={31},
       pages={1529\ndash 1551; arXiv:1512.02499},
}

\bib{profiles}{misc}{
      author={Hundertmark, F.},
       title={Profiles. {A}n algebraic approach to combinatorial connectivity},
         how={unpublished manuscript},
        date={2011},
}

\bib{GMX}{article}{
      author={Robertson, N.},
      author={Seymour, P.D.},
       title={Graph minors. {X}. {O}bstructions to tree-decomposition},
        date={1991},
     journal={J. Comb. Theory Ser. B},
      volume={52},
       pages={153\ndash 190},
}

\end{biblist}
\end{bibdiv}
